\documentclass[12pt,centertags,oneside]{amsart}
\usepackage{amsmath,amstext,amsthm,amscd,typearea,hyperref}
\usepackage{amssymb}
\usepackage{a4wide}
\usepackage[mathscr]{eucal}
\usepackage{mathrsfs}
\usepackage{typearea}
\usepackage{charter}
\usepackage{pdfsync}
\usepackage[a4paper,width=16.2cm,top=3cm,bottom=3cm]{geometry}

\numberwithin{equation}{section}

\newtheorem{theorem}{Theorem}[section]
\newtheorem{definition}[theorem]{Definition}
\newtheorem{proposition}[theorem]{Proposition}
\newtheorem{corollary}[theorem]{Corollary}
\newtheorem{lemma}[theorem]{Lemma}
\newtheorem{remark}[theorem]{Remark}

\newcommand{\cali}[1]{\mathscr{#1}}

\newcommand{\MA}{{\rm MA}}
\newcommand{\NMA}{{\rm NMA}}

\newcommand{\supp}{{\rm supp}}

\newcommand{\maxep}{\mathop{\mathrm{\rm  max_\epsilon}}\nolimits}

\newcommand{\const}{\mathop{\mathrm{const}}\nolimits}

\newcommand{\Leb}{\mathop{\mathrm{Leb}}\nolimits}
\newcommand{\Lip}{\mathop{\mathrm{Lip}}\nolimits}
\newcommand{\LLip}{{\mathop{\mathrm{\widetilde{Lip}}}\nolimits}}

\newcommand{\dist}{\mathop{\mathrm{dist}}\nolimits}
\renewcommand{\Re}{\mathop{\mathrm{Re}}\nolimits}

\newcommand{\vol}{\mathop{\mathrm{vol}}}

\newcommand{\eq}{\mathrm{eq} }

\newcommand{\ddc}{dd^c}
\newcommand{\dc}{d^c}

\def\mB{\mathcal{B}}
\def\mD{\mathcal{D}}
\def\mE{\mathcal{E}}

\def\mL{\mathcal{L}}

\def\mV{\mathcal{V}}
\def\mW{\mathcal{W}}

\newcommand{\PSH}{{\rm PSH}}

\newcommand{\ind}{{\bf 1}}

\newcommand{\exph}{{\rm exph}}

\newcommand{\Bc}{\cali{B}}
\newcommand{\Cc}{\cali{C}}

\renewcommand{\Mc}{\cali{M}}

\newcommand{\Pc}{\cali{P}}

\newcommand{\B}{\mathbb{B}}
\newcommand{\C}{\mathbb{C}}
\newcommand{\D}{\mathbb{D}}
\newcommand{\N}{\mathbb{N}}

\newcommand{\R}{\mathbb{R}}
\renewcommand\P{\mathbb{P}}

\renewcommand{\S}{\mathbb{S}}

\newcommand{\e}{\textbf{\textit{e}}}

\title[Equidistribution speed for Fekete points]{Equidistribution 
speed for Fekete points associated\\ with an ample line bundle}

\author{Tien-Cuong Dinh}
\address{Department of Mathematics, National University 
of Singapore, 10 Lower Kent Ridge Road, Singapore 119076.}
\email{matdtc@nus.edu.sg}
\thanks{T.-C.\ D.\  partially supported by Start-Up 
Grant R-146-000-204-133 from National University of
\break Singapore}
\author{Xiaonan Ma}
\address{Institut Universitaire de France
\&Universit\'e Paris Diderot - Paris 7,
UFR de Math\'ematiques, Case 7012,
75205 Paris Cedex 13, France.}
\email{xiaonan.ma@imj-prg.fr}
\thanks{X. M. partially supported by Institut Universitaire de France, 
ANR-14-CE25-0012-01 and funded through the Institutional Strategy 
of the University of Cologne within the German Excellence Initiative.}
\author{Vi{\^e}t-Anh Nguy{\^e}n}
\address{Math{\'e}matique-B{\^a}timent 425, UMR 8628, 
Universit{\'e} Paris-Sud, 91405 Orsay, France. 
}
\email{VietAnh.Nguyen@math.u-psud.fr}
\thanks{V.-A.\ N.\ partially supported by the Max-Planck institute 
for mathematics in Bonn.}
\date{June 11, 2015}
\begin{document}

\begin{abstract}
Let $K$ be the closure of a bounded open set with smooth
boundary in $\C^n$. 
A Fekete configuration of order $p$ for $K$ is a finite subset of 
$K$ maximizing the Vandermonde determinant associated with 
polynomials of degree $\leq p$.
 A recent theorem by Berman, Boucksom and  Witt Nystr\"om 
 implies that Fekete configurations for $K$ are asymptotically 
 equidistributed with respect to a canonical equilibrium measure, 
 as $p\to\infty$. We give here an explicit estimate 
 for the speed of convergence. The result also holds
 in a general setting of  Fekete points associated with  an ample 
  line bundle over a projective manifold.   Our approach requires 
a new estimate on  Bergman kernels for line bundles and 
quantitative results  in pluripotential  theory which are of
independent interest. 
\end{abstract}

\maketitle

\medskip

\noindent
{\bf Classification AMS 2010}: 32U15 (32L05).

\medskip

\noindent
{\bf Keywords:} Fekete points, equilibrium measure,  equidistribution, Bergman kernel,\break
Monge-Amp\`ere operator, Bernstein-Markov property.

\bigskip
{\small 
\noindent
{\sc R\'esum\'e.} Soit $K$ l'adh\'erence d'un ouvert born\'e \`a 
bord lisse dans $\C^n$. 
Une configuration de  Fekete d'ordre $p$ pour $K$ est 
un sous-ensemble fini de 
$K$ qui maximise le d\'eterminant de Vandermonde associ\'e aux 
polyn\^omes de degr\'e  $\leq p$.
Un th\'eor\`eme r\'ecent de  Berman, Boucksom et  Witt Nystr\"om
implique que les configurations de 
Fekete sont asymptotiquement \'equir\' eparties par rapport \`a 
une mesure d'\'equilibre canonique quand  $p\to\infty$. 
Nous donnons ici une estimation pr\'ecise de la vitesse 
de convergence. Le r\'esultat est aussi valable dans un cadre 
g\'en\'eral des points de   Fekete associ\'es \`a un fibr\'e en 
droites ample au-dessus d'une vari\'et\'e projective. 
Notre approche n\'ecessite une estimation nouvelle sur 
les noyaux de   Bergman pour les fibr\'es en droites et 
des r\'esultats quantitatifs de la th\'eorie du pluripotentiel 
qui sont d'int\'er\^et ind\'ependant.

\medskip\noindent
{\bf Classification de l'AMS 2010} : 32U15 (32L05).

\medskip

\noindent
{\bf Mots-cl\'es:} points de Fekete, mesure d'\'equilibre,  
\'equidistribution, noyau de Bergman, op\'erateur de
Monge-Amp\`ere, propri\'et\'e de Bernstein-Markov.

}

\tableofcontents


\addcontentsline{toc}{section}{Notation and some remarks}

\noindent
{\bf Notation.} Throughout the paper, $L$ denotes an ample  
holomorphic line bundle over a projective manifold $X$ 
of dimension $n$. Fix also a smooth Hermitian metric $h_0$ 
on $L$ whose first Chern form, 
denoted by $\omega_0$, is a K\"ahler form. For simplicity, 
we use the K\"ahler metric on $X$ induced by $\omega_0$. 
The induced distance is denoted by $\dist$.   
Define
$\mu^0:=\|\omega_0^n\|^{-1}\omega_0^n$ the probability measure 
associated with the volume form $\omega_0^n$. 
The space of holomorphic sections of $L^p:=L^{\otimes p}$,
the $p$-th power of $L$, is denoted by $H^0(X,L^{p})$.
Its dimension is denoted by $N_p$. 
The metric $h_0$ induces, in a canonical way, metrics 
on the line bundle $L^{p}$ over $X$, 
the vector bundle of the product 
$L^{p}\times \cdots \times L^{p}$ ($N_p$ times)
 over $X^{N_p},$ and the determinant of the last one which is 
 a line bundle over $X^{N_p}$ and denoted by 
 $(L^p)^{\boxtimes N_p}$. 
 For simplicity, the norm, induced by $h_0$, of a section 
 of these vector bundles is denoted by $|\cdot |$. 

A general singular metric on $L$ has the form $h=e^{-2\psi}h_0$, 
where $\psi$ is an integrable function on $X$  with values in 
$\R\cup\{\pm\infty\}$. Such a  function $\psi$ is  called  
a {\it weight}. 
 It also induces singular metrics on the above vector bundles,
and we denote by $|\cdot|_{p\psi}$ the corresponding norm of 
a section 
of $L^{p}$ or the associated determinant line bundle over 
$X^{N_p}$. This is a function on $X$ or $X^{N_p}$ respectively.
If $K$ is a subset of $X$, the supremum on $K$ 
or $K^{N_p}$ of this function is denoted by 
$\|\cdot\|_{L^\infty(K,p\psi)}$ or 
$\|\cdot\|_{L^\infty(K^{N_p},p\psi)}$. 
Its $L^2(\mu)$ or $L^2(\mu^{\otimes N_p})$-norm  
is denoted by $\|\cdot\|_{L^2(\mu,p\psi)}$ or 
$\|\cdot\|_{L^2(\mu^{\otimes N_p},p\psi)}$, where $\mu$ is a
probability measure on $X$.  We sometimes drop the power $N_p$
for simplicity. 
In the same way, we often add the index ``$\psi$" or ``$p\psi$",
if necessary, to inform the use of the weight $\psi$
 for $L$ and hence $p\psi$ for $L^{p}$. 

The notations $\rho_p(\mu,\phi)$, $\Bc_p(\mu,\phi)$ will 
be introduced in Subsection \ref{section_Bergman}, 
$\mB^\infty_p(K,\phi)$,\break 
$\mB^2_p(\mu,\phi)$,  $\mL_p(K,\phi)$, $\mL_p(\mu,\phi)$, 
$\mE(\phi)$, $\mE_\eq(K,\phi)$
 in Subsection \ref{subsection_preliminaries}, and 
$\mV_p(\phi_1,\phi_2)$, $\mW(\phi_1,\phi_2)$, $\epsilon_p$,
$\mD_p(K,\phi)$ 
 in Subsection \ref{subsection_Main_estimates}. 
 Let $\B(x,r)$ denote the ball of center $x$ and radius $r$ in $X$ 
 or in an Euclidean space. Similarly, $\D(x,r)$ is the disc of center 
 $x$ and radius $r$ in $\C$, $\D_r:=\D(0,r)$ and $\D:=\D(0,1)$.  
 The Lebesgue measure on an Euclidean space is denoted 
 by $\Leb$. 
The operators $\dc$ and $\ddc$ are defined by 
$$\dc:= {\sqrt{-1}\over 2\pi }(\overline\partial -\partial) 
\quad \text{and} \quad
\ddc:={\sqrt{-1}\over \pi} \partial\overline\partial.$$
For $m\in\N$ and $0<\alpha\leq 1$,  $\Cc^{m,\alpha}$ is
the class of $\Cc^m$ functions/differential forms whose 
partial derivatives of 
order $m$ are H\"older continuous with H\"older exponent 
$\alpha$. We have $\Cc^{m,\alpha}=\Cc^{m+\alpha}$ except for 
$\alpha=1$.  We use the natural norms on these spaces and 
for simplicity, define $\|\cdot\|_m:=1+\|\cdot\|_{\Cc^m}$ and 
$\|\cdot\|_{m,\alpha}:=1+\|\cdot\|_{\Cc^{m,\alpha}}$.  
 Denote by $\Lip$ the space of Lipschitz functions which is also
 equal to $\Cc^{0,1}$ and by $\LLip$ the space of functions $v$ 
 such that
$|v(x)-v(y)|\lesssim - \dist(x,y)\log\dist (x,y)$ 
for $x,y$ close enough. We endow the last space with the norm
$$\|v\|_{\LLip}:=\|v\|_\infty+\inf \big \{A\geq 0: 
\quad  |v(x)-v(y)|\leq - A\dist(x,y)\log\dist (x,y) \text{ if } 
\dist(x,y)\leq 1/2\big\}.$$

 A function $\phi:X\to \R\cup\{-\infty\}$ is  called 
{\it quasi-plurisubharmonic} (quasi-p.s.h. for short) if it is locally the  
sum of a plurisubharmonic (p.s.h. for short) and  a smooth function. 
A quasi-p.s.h.  function $\phi$ is called {\it $\omega_0$-p.s.h.}
 if $\ddc\phi+\omega_0\geq 0$ in the sense of currents. 
 Denote by $\PSH(X,\omega_0)$ the set of such functions. 
 If $\phi$ is a bounded function in $\PSH(X,\omega_0)$, define
  {\it the associated Monge-Amp\`ere measure} 
  and {\it normalized Monge-Amp\`ere measure}  by
$$\MA(\phi):=(\ddc\phi+\omega_0)^n \quad \text{and}
 \quad \NMA(\phi):=\|\MA(\phi)\|^{-1}\MA(\phi).$$
So $\MA(\phi)$ is a positive measure and $\NMA(\phi)$ is 
a probability measure on $X$.  
A quasi-p.s.h.  function $\phi$ is called 
{\it strictly $\omega_0$-p.s.h.} if $\ddc\phi+\omega_0$ 
is larger than a K\"ahler form in the sense of currents, 
see \cite{Demailly12, DinhSibony4} 
for the basic notions and results of pluripotential theory. 

\medskip
\noindent
{\bf Some remarks.} The constants involved in our computations 
below may depend on $X,L,h_0$ and hence on $\omega_0$ and 
$\mu^0$. However, they do not depend on the other weights
used for the line bundle $L$ 
but only on the upper bounds of  suitable norms 
($\Cc^\alpha$, $\LLip$, ...) of these weights. 
This property can be directly seen in our arguments.
For simplicity, we will not repeat it in each step of the proofs. 
The notations $\gtrsim$ and $\lesssim$ mean inequalities up to a 
positive multiple constant.

\section{Introduction} \label{introduction}

Let $K$ be a non-pluripolar compact subset of $\C^n$. 
{\it The pluricomplex Green function} of $K$, denoted by 
$V_K^*(z)$, 
  is the upper-semicontinuous regularization of the
  Siciak-Zahariuta extremal function
$$V_K(z):=\sup \big\{ u(z):  u \text{ p.s.h. on }\C^n,  
u|_{K}\leq 0, u(w)-\log\|w\|=O(1) \text{ as } w\to\infty\big\}.$$
This function $V^*_K$ is locally bounded, p.s.h.  and 
$(\ddc V_K^*)^n$ defines a probability measure with support in $K$.
 It is called {\it the equilibrium measure} of $K$ and denoted by 
 $\mu_\eq(K)$, see \cite{Siciak81,Zaharjuta}.

Let $\Pc_p$ be the set of holomorphic polynomials of degree 
$\leq p$ on $\C^n$. This is a complex vector space of dimension 
$$N_p:=\binom{p+n}{n} ={1\over n!}p^n +O(p^{n-1}).$$
Let $(e_1,\ldots, e_{N_p})$ 
be a basis of $\Pc_p$. 
Define for $P=(x_1,\ldots,x_{N_p})\in (\C^n)^{N_p}$
 {\it the Vandermonde determinant} $W(P)$ by 
$$W(P):=\det\left(\begin{array}{ccc}
e_1(x_1) & \ldots & e_1(x_{N_p}) \\
\vdots & \ddots & \vdots \\
e_{N_p}(x_1) & \ldots & e_{N_p}(x_{N_p})
\end{array} \right).
$$
A point $P\in K^{N_p}$ is called {\it a Fekete configuration for $K$}
if the function $|W(\cdot)|$,
 restricted to $K^{N_p}$, achieves its maximal value at $P$. 
 It is not difficult to check that this definition does not depend 
 on the choice of 
 the basis $(e_1,\ldots, e_{N_p})$, see  \cite{SaffTotik}.

Recently,  Berman, Boucksom and  
Witt Nystr\"om have  proved that Fekete points\break 
$x_1,\ldots,x_{N_p}$ are asymptotically 
equidistributed with respect to the equilibrium measure 
$\mu_\eq(K)$ as $p$ tends to infinity \cite{BBW}. 
This property had been conjectured for  quite some time, 
probably going back to the pioneering  work of  
Leja in \cite{Leja57a, Leja57b}, where the dimension 1 case 
was obtained. See also \cite{Levenberg, SaffTotik}  
for more recent references on this topic. More precisely, let
$$\mu_p:={1\over N_p} \sum_{j=1}^{N_p} \delta_{x_j}$$
denote the probability measure equidistributed on 
$x_1,\ldots, x_{N_p}$. We call it 
{\it a Fekete measure of order $p$}. 
The above equidistribution result says that in the weak-$*$ topology
$$\lim_{p\to\infty} \mu_p=\mu_\eq(K).$$
In fact, this theorem by Berman, Boucksom and Witt Nyst\"om holds 
in a more  general  context of Fekete points associated with 
a line bundle. We will discuss this case later together with 
an interesting new approach by Ameur, Lev and Ortega-Cerd\`a 
\cite{AmeurOrtega, LevOrtega}.  

Fekete points are well known to be useful in several problems
in mathematics and mathematical physics. It is therefore important 
to study  the speed of the above convergence. For this purpose, 
it is necessary to make some hypothesis on the compact set $K$. 
For instance, we have the following result, 
see also Corollary \ref{cor_main_K}.

\begin{theorem} \label{th_main_Cn}
Let $K$ be the closure of a bounded non-empty open subset of 
$\C^n$ with $\Cc^2$ boundary.
Then for all  $0<\gamma\leq 2$ and $\epsilon>0$, there is 
a constant $c=c(K,\gamma,\epsilon)>0$, independent of $p>1$, 
such that 
$$|\langle \mu_p-\mu_\eq(K),v\rangle|
\leq c \|v\|_{\Cc^\gamma} p^{-\gamma/36+\epsilon} $$
for every Fekete measure $\mu_p$ of order $p$ and 
every test function $v$ of class $\Cc^\gamma$ on $\C^n$.
\end{theorem}

In fact, our result is still true 
in  a more general setting that we will state below after introducing 
necessary notation and  terminology.

\medskip

Let $L$ be an ample holomorphic line bundle over 
a projective manifold $X$ 
of dimension $n$. Fix a smooth Hermitian metric $h_0$ on $L$ 
whose first Chern 
form $\omega_0:= \frac{\sqrt{-1}}{2\pi}R^L_{0}$ is a K\"ahler form,
where $R^L_{0}$ is the curvature of the Chern connection on 
$(L,h_{0})$. 

\begin{definition} \rm \label{def_weighted_K}
We call {\it weighted compact subset of $X$} a data $(K,\phi)$, 
where $K$ is a non-pluripolar compact subset of $X$ and 
$\phi$ is a real-valued continuous function on $K$. 
The function $\phi$ is called {\it a weight} on $K$. 
The {\it equilibrium weight} associated with $(K,\phi)$ is
 the upper semi-continuous regularization $\phi_K^*$ of 
 the function 
$$ \phi_K(z):=\sup\big\{\psi(z) : \psi \ \  \omega_0 \text{-p.s.h.}, 
\ \psi\leq \phi \text{ on } K\big\}.$$
We call {\it equilibrium measure of $(K,\phi)$} 
the normalized Monge-Amp\`ere measure 
$$\mu_\eq(K,\phi):=\NMA(\phi_K^*).$$
\end{definition}

Note that the equilibrium measure $\mu_\eq(K,\phi)$ 
is a probability measure
 supported by $K$ and $\phi_K^*=\phi_K$ almost everywhere 
 with respect to this measure, see e.g., \cite{BermanBoucksom}. 

\begin{definition} \rm \label{def_K_regular}
Denote by $P_K$ the projection onto $\PSH(X,\omega_0)$ 
which associates $\phi$ with $\phi_K^*$. We say that $(K,\phi)$
 is {\it regular} if $\phi_K$ is upper semi-continuous,
 i.e., $P_K\phi=\phi_K$. Let $(E, \|\ \|_E)$ be  a normed vector 
 space of functions on $K$ and $(F,\| \ \|_F)$ a normed vector 
 space of functions on $X$. We say that $K$ is {\it $(E,F)$-regular}  
if $(K,\phi)$ is regular for  $\phi\in E$ 
and if the  projection $P_K$  sends bounded subsets of $E$ into 
bounded subsets of $F$. 
\end{definition}

We will see in Theorem \ref{th_reg_K} below that when $K$ is 
the closure of an open set with $\Cc^2$ boundary, then it is
 $(\Cc^\alpha,\Cc^\alpha)$-regular for $0<\alpha<1$, 
 i.e., $(E,F)$-regular with $E=\Cc^\alpha(K)$ and $F=\Cc^\alpha(X)$. 

Consider now an integrable real-valued function $\psi$ on $X$
and the singular Hermitian metric
$h:=e^{-2\psi}h_0$ on the line bundle $L$. 
We will use the notations given at the beginning of the paper. 
Consider also a basis $S_p=(s_1,\ldots,s_{N_p})$ of the vector 
space $H^0(X,L^{p})$, where $N_p:=\dim H^0(X,L^{p})$. 
This basis can be seen as a section of the rank $N_p$ 
vector bundle of the product $L^{p}\times \cdots\times L^{p}$ 
 ($N_p$-times) 
 over $X^{N_p}$. The determinant line bundle associated with 
 this vector bundle is denoted by $(L^p)^{\boxtimes N_p}$.
The determinant $\det(s_i(x_j))_{1\leq i,j\leq N_p}$ for 
$P=(x_1,\ldots,x_{N_p})$ in $X^{N_p}$ defines a section of 
the last line bundle over $X^{N_p}$ that we will denote by
$\det S_p$ or $\det(s_i(x_j))$. The metric $h_0$ induces in 
a canonical way a metric $(h_0^p)^{\boxtimes N_p}$ on  
$(L^p)^{\boxtimes N_p}$.
As mentioned above, we denote by $|\det(s_i(x_j))|$ the norm of 
$\det(s_i(x_j))$ with respect to $(h_0^p)^{\boxtimes N_p}$.
For $P=(x_1,\ldots,x_{N_p})$ in $X^{N_p}$, we will consider 
{\it the weighted Vandermonde determinant} 
$$|\det(s_i(x_j))|_{p\psi}:=|\det(s_i(x_j))|e^{-p\psi(x_1)
-\cdots-p\psi(x_{N_p})}.$$
The following notion does not depend on the choice of
the basis $S_p=(s_1,\ldots,s_{N_p})$. 

\begin{definition} \rm \label{def_Fekete_L}
The point $P=(x_1,\ldots,x_{N_p})$ in $K^{N_p}$ is called 
{\it a Fekete configuration of order $p$ of $(L,h_0)$} in the 
weighted compact set $(K,\phi)$ if the above  weighted 
Vandermonde determinant, restricted to $K^{N_p}$, 
achieves its maximal value at $P$. 
The associated probability measure 
$${1\over {N_p}}(\delta_{x_1}+\cdots +\delta_{x_{N_p}}),$$
on $K$ is called {\it a Fekete measure of order $p$}. 
\end{definition}

In order to study the speed of equidistribution of Fekete points,
 it is convenient to use some distance notions on the space 
 $\Mc(X)$ of (Borel) probability measures on $X$. 
For $\gamma>0$, define
 the distance $\dist_\gamma$ between two measures 
   $\mu$ and $\mu'$ in $\Mc(X)$ by 
$$\dist_\gamma(\mu,\mu'):=\sup_{\| v\|_{\Cc^\gamma}\leq 1} 
\big |\langle \mu-\mu', v\rangle\big|,$$
where $v$ is a test smooth real-valued function. 
This distance induces the weak topology on $\Mc(X)$. 
By interpolation between Banach spaces 
(see \cite{DinhSibony4, Triebel}), 
for $0<\gamma\leq \gamma'$, there exists $c>0$ such that
\begin{align}\label{eq:n1.23}
\dist_{\gamma'}\leq \dist_\gamma
\leq c[\dist_{\gamma'}]^{\gamma/\gamma'}.
\end{align}
Note that $\dist_1$ is equivalent to the classical 
Kantorovich-Wasserstein distance.

Here is our main result which is the version of 
Theorem \ref{th_main_Cn} in the general setting.
It is already interesting for $K=X$.

\begin{theorem}\label{thm_main}
Let $X,L,h_0$ be as above and $K$ a non-pluripolar compact 
subset of $X$.
Let $0<\alpha\leq 2,$ $ 0<\alpha'\leq 1$ and $0<\gamma\leq 2$ 
be constants. 
Assume that $K$ is $(\Cc^\alpha,\Cc^{\alpha'})$-regular.
 Let $\phi$ be a $\Cc^\alpha$ real-valued function on $K$ 
 and $\mu_\eq(K,\phi)$ the equilibrium measure associated 
 with the weighted set $(K,\phi)$. Then, there is $c>0$ 
 such that for every $p>1$ and every Fekete measure $\mu_p$ 
 of order $p$ associated with $(K,\phi)$, we have 
$$\dist_\gamma(\mu_p, \mu_\eq(K,\phi))
\leq c p^{-\beta\gamma}(\log p)^{3\beta\gamma}\quad \text{with} 
\quad\beta:=  \alpha'/(24+12\alpha').$$
 \end{theorem}

We will see later in Theorem \ref{th_reg_K} that the hypothesis on 
$K$ is satisfied for $\alpha=\alpha'<1$ when $K$ is the closure of 
an open set with $\Cc^2$ boundary (we think that the techniques
we use can be applied to study other classes of compact sets but 
we don't develop this direction here). So the result  below is 
a consequence of Theorem \ref{thm_main} for $\alpha=\alpha'<1$. 

\begin{corollary} \label{cor_main_K}
Let $X,L,h_0$ be as above and $K$ the closure of a non-empty 
open subset of $X$ with $\Cc^2$ boundary. Let $\phi$  
be a $\Cc^\alpha$ real-valued function on $K$, $0<\alpha<1$, and  
$\mu_\eq(K,\phi)$ the equilibrium measure associated 
 with $(K,\phi)$. Then, for every $0<\gamma\leq 2$, there is $c>0$ 
 such that for every $p>1$ and every Fekete measure $\mu_p$ 
 of order $p$ associated with $(K,\phi)$, we have 
$$\dist_\gamma(\mu_p, \mu_\eq(K,\phi))
\leq c p^{-\beta\gamma}(\log p)^{3\beta\gamma}\quad \text{with} 
\quad\beta:=  \alpha/(24+12\alpha).$$
\end{corollary}

When $X$ is the projective space $\P^n$ and $L$ is the tautological
 line bundle $\mathcal O(1)$ on $\P^n$, we can consider $X$ as 
 the natural compactification of $\C^n$ and the sections in 
 $H^0(X,L^{p})=H^0(\P^n,\mathcal O(p))$ can be identified to 
 polynomials of degree $\leq p$ on $\C^n$. We then see that 
 Theorem \ref{th_main_Cn} is a particular case of the last corollary.

Our theorem applies to the case where $K=X$ and $\phi$ is 
a smooth function on $X$. If the metric $h:=e^{-2\phi}h_0$ of $L$ 
has strictly positive curvature form, our approach gives an estimate
 better than the one in the last theorem. Namely, we have the
 following result, see also Remark \ref{rk_main_X_positive}.

\begin{theorem}\label{thm_main_X_positive}
Let $X,L$ and $h_0$ be as above. Let $\phi$ be a $\Cc^{3}$ 
real-valued function on $X$ such that the first Chern form of 
the metric $h:=e^{-2\phi}h_0$ is strictly positive. 
Let $\mu_\eq(X,\phi)$ denote the equilibrium measure
associated with the weighted set $(X,\phi)$. 
Then for any $0<\gamma\leq 3,$ there is $c>0$ 
such that
$$\dist_\gamma(\mu_p,\mu_\eq(X,\phi)) 
\leq c p^{-\gamma/12}(\log p)^{\gamma/4}$$
for all $p>1$ and all Fekete measures $\mu_p$ of order  
$p$ associated with $(X,\phi)$.
\end{theorem}

This result is close to the one recently obtained by  
Lev and Ortega-Cerd\`a  in \cite{LevOrtega}.  
These authors proved that when $\phi$ is  
smooth $\omega_0$-strictly p.s.h., there is
a constant $c>0$ such that
\begin{align}\label{eq:n1.27}
c^{-1} p^{-1/2}\leq \dist_1(\mu_p, \mu_\eq(X,\phi))
\leq c p^{-1/2}  
\end{align}
for all $p$ and  Fekete measures $\mu_p$ of order 
$p$ associated with $(X,\phi)$.
Using \eqref{eq:n1.23}, 
we can deduce  similar estimates for $\dist_\gamma$ with 
$0<\gamma\leq 1$.
 So the result of Lev and Ortega-Cerd\`a is  optimal for 
 $0<\gamma\leq 1$ in their assumption.
Although  for $0<\gamma\leq 1$ estimate in
Theorem \ref{thm_main_X_positive}.
is weaker than \eqref{eq:n1.27} and its interpolated version,  
our assumption  of smoothness for $\phi$ is only $\Cc^{3}$ 
and can be easily reduced to $\Cc^\alpha$  with  similar estimates 
depending on $\alpha$, see Remark \ref{rk_main_X_positive}.
Of course, in the case where the curvature of the metric induced 
by $\phi$ is only semi-positive or even not semi-positive, 
one can apply Corollary \ref{cor_main_K} to $K=X$. 

\smallskip

In their  approach, Lev and Ortega-Cerd\`a  relate the  
equidistribution of Fekete points to the problem of sampling 
and interpolation on line  bundles as in a previous work by  
 Ameur and Ortega-Cerd\`a \cite{AmeurOrtega}. 
The main ingredients of their method consist in using  
Toeplitz operators as well as  known asymptotic expansions for 
the Bergman kernels on/off the diagonal  of $X\times X$ due to
\cite{Catlin,Lindholm,Tian,Zelditch}, 
  cf. also \cite{MaMarinescu07, MaMarinescu14}.
  The key points here are  (1)  the Fekete configurations
  are also  sampling and  interpolation,  and (2)
the points of such a configuration are geometrically  equidistributed.
These  crucial properties  are obtained using the assumption that
the  metric weight $\phi$ is  smooth $\omega_0$-strictly p.s.h.  
  
 \smallskip

Our approach is different because our metric weight $\phi$ is, 
in general, only H\"older continuous and  it may originally be 
defined on a  proper compact set $K\subset X.$ In this context, 
$P_K\phi$ is only weakly $\omega_0$-p.s.h., and
moreover,  not smooth in general.
So  the result by Lev and Ortega-Cerd\`a  is not applicable  
in the general context.

We will follow the original method of
Berman, Boucksom and Witt Nystr\"om \cite{BermanBoucksom, BBW}.
We will need, among other things, a controlled regularization for 
quasi-p.s.h. functions, quantitative properties of  quasi-p.s.h. 
envelopes of functions and an estimate of  
Bergman kernels associated with holomorphic line bundles. 
These results are of independent interest and will be presented 
in the next section while
the proofs of the main results will be given in the last section.

 \medskip\noindent
{\bf Acknowledgment.} 
The paper was partially written during the visits of the second and  
the third authors at National University of Singapore, University of 
Cologne and Max-Planck institute for mathematics in Bonn. 
They would like to thank these organizations for their very warm 
hospitality.

\section{Quasi-p.s.h. functions, equilibrium weight and 
Bergman functions}  \label{section_envelope}

Let $X$ be a compact  K\"ahler manifold of dimension $n$ and
 let $\omega_0$ be a fixed K\"ahler form on $X$. 
We will use later  the equilibrium weight $P_K\phi$
associated with a regular weighted compact set $(K,\phi)$ of $X$. 
This is a quasi-p.s.h. function which is not smooth in general.
 So we will need to approximate it by smooth quasi-p.s.h. functions 
 and control the cost of this regularization procedure. 

In this section, we will give a version of the theorem of regularization
 for H\"older continuous quasi-p.s.h. functions and study 
 the H\"older continuity of equilibrium weights.
The behavior of Bergman functions associated with the powers of 
a line bundle with small positive curvature is crucial in our approach. 
This question will also be considered here in the last subsection. 

\subsection{Regularization of quasi-p.s.h. functions}
\label{subsection_regularization}

The purpose of this subsection is to  establish the following
 regularization theorem for H\"older continuous quasi-p.s.h. functions  
 with a control of positivity and controlled $\Cc^m$ norms.
 
\begin{theorem}\label{thm_Demailly}
For each $0<\alpha\leq 1,$ 
there exist $c> 0$ which only depends  on $X,$ $\omega_0,\alpha,$ and
$c_m>0$   
 which only depends  on $X,\omega_0,\alpha$ and $ m\in\N^*$  
 satisfying the following property.
Let $\phi$ be  an  $\omega_0$-p.s.h. function on $X$ of class 
$\Cc^{0,\alpha}.$  Then,  for each $0<\epsilon\leq 1,$ there  exists 
a smooth function $\phi_\epsilon$  such that
 \begin{enumerate}
 \item[a)] $\phi_\epsilon$ is $\omega_0$-p.s.h.; 
 \item[b)]  $\|\phi_\epsilon -\phi\|_\infty
 \leq c \epsilon^{\alpha} \|\phi\|_{0,\alpha}$
 (see the beginning of the paper for notation);
 \item[c)] $\| \phi_\epsilon\|_{\Cc^m(X)}
 \leq c_m \epsilon^{-m+\alpha} \|\phi\|_{0,\alpha}$ for $m\in \N^*.$
 \end{enumerate}
\end{theorem}

We  are   inspired  by  Demailly's regularization theorem
\cite{Demailly12, Demailly94} and  a technique of
Blocki-Kolodziej  \cite{BlockiKolodziej}.
 First,  we construct  suitable regularized  maximum functions. 
 Fix a function $\vartheta\in\Cc^\infty(\R,\R^+)$ 
 with support in $[-1,1]$
 such that $\int_\R \vartheta(h)dh=1$ and  
 $\int_\R h\vartheta(h)dh=0.$ 
 For each  $0<\epsilon\leq 1$ and each  integer $l\geq 1,$ 
consider the {\it regularized  maximum  function}  
$\max_\epsilon:\ \R^l \to \R$ defined by
$$ \maxep (t_1,\ldots,t_l):=\int_{\R^l}  
\max(t_1+h_1,\ldots,t_l+h_l)\epsilon^{-l}
 \prod_{i=1}^l\vartheta(h_i/\epsilon)dh_1\ldots dh_l. $$
Here are 
some  properties of  $\max_\epsilon$  which will be used later.
The notation $(t_1,\ldots, \widehat{t_i},\ldots,t_l)$ below 
means that the component $t_i$ is omitted in the expression.

\begin{lemma}\label{lem_properties_regularized_max}  
\begin{enumerate}
\item[a)]   $\max_\epsilon (t_1,\ldots,t_l)$ is  
non-decreasing in all variables, 
smooth and convex  on $\R^l;$
\item [b)] $\max (t_1,\ldots,t_l)\leq \max_\epsilon (t_1,\ldots,t_l)
\leq \epsilon+\max (t_1,\ldots,t_l);$  
\item[c)]  $\max_\epsilon (t_1,\ldots,t_l) 
=\max_\epsilon (t_1,\ldots, \widehat{t_i},\ldots,t_l)$ 
if $t_i+2\epsilon \leq \max(t_1,\ldots, \widehat{t_i},\ldots,t_l);$
\item[d)] if $u_1,\ldots,u_l$ are p.s.h. functions defined on some  
domain $D$ in $\C^n,$ then so is $\max_\epsilon (u_1,\ldots,u_l).$   
\item[e)] If $u_1,\ldots,u_l$ are real-valued  functions in 
$\Cc^m(D),$  where   
$m\in\N^*$ and $D$ is a domain in $\C^n,$  then there is 
a constant $c_{l,m}>0$ depending only on $l, m$ and $\vartheta$ 
such that
$$\left\|  \maxep (u_1,\ldots,u_l) \right\|_{\Cc^m}\leq  \epsilon
+ \sup_{1\leq i\leq l} \| u_i\|_\infty +c_{l,m}  \sum_{r_{ij}}
  \epsilon^{1-\sum r_{ij}} \prod_{i,j}     \| u_i\|_{\Cc^j}^{r_{ij}} ,$$
the sum being taken over all $r_{ij}>0$ with $1\leq i\leq l$ and
$j\geq 1$ such that  $\sum j r_{ij}\leq  m.$  
\end{enumerate} 
 \end{lemma}
\proof 
Assertions a)-d)  are contained in Lemma I.5.18 of \cite{Demailly12}, 
where the above properties of $\vartheta$ are used.
We turn to assertion e). 
Note that assertion b) allows us to bound the sup-norm of 
$ \maxep (u_1,\ldots,u_l)$, and hence explains the presence of 
$\epsilon + \sup\| u_i\|_\infty$ in assertion e).

Observe that the function $\max$ is Lipschitz. Therefore, 
any partial derivative of order 1 of $ \maxep (u_1,\ldots,u_l)$, 
seen as a function in $D$, is a finite sum of integrals of type 
\begin{equation} \label{eq:max:der1}
v\int_{\R^l} \Phi(u_1+h_1,\ldots,u_l+h_l)\epsilon^{-l} 
\prod_{i=1}^l\vartheta(h_i/\epsilon)dh_1\ldots dh_l, 
\end{equation}
where $\Phi$ is a partial derivative of order 1 of $\max$ and $v$ 
is a partial derivative of order 1 of a function $u_i$. 
Note that $\Phi$ is bounded. 

Performing the change of variables 
$u_i+h_i=s_i$, the expression in \eqref{eq:max:der1} is equal to
$$v\int_{|s_i-u_i|\leq\epsilon} \Phi(s_1,\ldots,s_l) 
\epsilon^{-l}\prod_{i=1}^l\vartheta\big({s_i-u_i\over \epsilon}
\big)ds_1\ldots ds_l,$$
which is a function in $D$. 
We see that any derivative up to order $m-1$ of this function
is bounded by a constant times
$$\sum_{r_{ij}}  \epsilon^{1-\sum r_{ij}} \prod_{i,j}   
\| u_i\|_{\Cc^j}^{r_{ij}},$$
where the sum is taken over all $r_{ij}>0$ with $1\leq i\leq l$ 
and $j\geq 1$ such that $\sum jr_{ij}\leq m$.  
This, together with the control of the sup-norm using b), implies 
assertion e).
 \endproof

Recall the following standard  regularization  by convolution.
Let $\rho(z):=\hat\rho(|z|)\in\Cc^\infty_0(\C^n)$ be  
a radial function such that 
$\hat\rho\geq 0,$ 
 $\hat\rho(t)=0$ for $t\geq 1,$ $\int_{\C^n}\rho  d\Leb=1,$
 where   $\Leb$ is the  Lebesgue measure on $\C^n.$
For $\delta>0$  we set  
$\rho_\delta(z):=\delta^{-2n}\rho(z/\delta).$ 
For every function $u$ on an open set $U\subset \C^n$  and 
every  subset  $U'\Subset U,$  define
\begin{align}\label{eq:n3.7}
u_\delta(z):=(u*\rho_\delta)(z)
=\int_{\C^n}  u(z-\delta w)\rho(w)d\Leb (w) \quad \text{with} 
\quad z\in U',
\end{align}
for $0<\delta <\dist(U',b U).$
If  $u$ is in $\Cc^{0,\alpha}(U)$  then  $u_\delta$ is in 
$\Cc^\infty(U')$  
and  we have  
\begin{align}\label{eq:n3.9}
\|u_\delta  -u\|_{\infty, U'}\lesssim \|u\|_{\Cc^{0,\alpha}}
\delta^{\alpha}\quad\text{and}
\quad  
\|u_\delta  \|_{\Cc^m(U')}\lesssim  \|u\|_{\Cc^{0,\alpha}}
\delta^{-m+\alpha} \quad \text{for} \quad m\in\N^*.
\end{align}
If $u$ is p.s.h. then $u_\delta$ is  also  p.s.h. and $u_\delta$ 
is decreasing to $u$ 
as $\delta\searrow 0.$  We need the following elementary lemma, 
whose proof is left to the reader, see also \cite{BlockiKolodziej}.

\begin{lemma}\label{lem_regularization}
Let  $F:\ W\to W'$ be a  biholomorphic map between two open 
subsets $W$ and $W'$ of $\C^n$.
Let  $u\in \PSH(W)\cap\Cc^{0,\alpha}(W)$ with $0<\alpha\leq 1$.  
Then, for every set $U\Subset W$ we can find a constant 
$\delta_U>0$  such that for $0<\delta  <\delta_U,$
the function
$u^F_\delta:=(u\circ F^{-1})_\delta\circ F$ is  well-defined  on  
a neighborhood  of $\overline U.$
Moreover, there  are  $c_U>0$ and $c_{U,m}>0$ for $m\in\N^*$ 
such that when $0<\delta<\delta_U,$  
$$\|u^F_\delta  -u\|_{\infty,U}\leq  c_U\|u\|_{\Cc^{0,\alpha}}
\delta^{\alpha} \quad\text{and}\quad  
\|u^F_\delta \|_{\Cc^m(U)}\leq c_{U,m}\|u\|_{\Cc^{0,\alpha}}
\delta^{-m+\alpha}.$$
\end{lemma}

\medskip

\noindent {\bf  End of the proof of Theorem  \ref{thm_Demailly}.}
Denote for simplicity $M:=\|\phi\|_{0,\alpha}$. 
The constants we will use below do not depend on $M$. 
Observe that  we only need to construct a 
$(1+c'M\epsilon^\alpha)\omega_0$-p.s.h. function $\phi_\epsilon$ 
such that 
\begin{align}\label{eq:n3.12}
\|\phi_\epsilon-\phi\|_\infty
\leq cM\epsilon^\alpha
 \quad \text{and} \quad \|\phi_\epsilon\|_{\Cc^m}
  \leq c_m M\epsilon^{-m+\alpha}
  \text{ for } m\geq 1,
\end{align}
where $c,c'$ and $c_m$ are constants. 
Indeed, we can just multiply it by 
$(1+c'M\epsilon^\alpha)^{-1}$
 in order to obtain a function as in Theorem  \ref{thm_Demailly}.
We can also add to this function a constant times
$M\epsilon^\alpha$ if we want to get a function larger or smaller 
than $\phi$.

First fix a finite cover of $X$ by small enough local charts 
$(U_j)_{j\in J}.$ We also choose a  
 finite cover of $X$ by local charts $(V_j)_{j\in J}$ indexing  
by the same index set $J$ such that $V_j\Subset U_j$.  
For each $j\in J$ fix a smooth function $f_j$ 
defined on   a neighborhood of $\overline U_j$ such that
\begin{align}\label{eq:n3.15}
  \ddc f_j= \omega_0   
\quad\text{on a neighborhood of}\quad \overline U_j.   
\end{align}
Then  the function  
 \begin{align}\label{eq:n3.16}  u_j:=\phi+f_j
 \end{align}
satisfies  $
\ddc u_j=\ddc \phi+\ddc f_j=\ddc\phi+\omega_0
\geq   0.$
So $u_j$ is p.s.h. on $U_j.$

Let $j$ and $k$ be in $J$ such that $U_j\cap U_k\not=\varnothing.$
There are two natural ways  to  regularize  
the restriction $u_j|_{U_j\cap U_k}$ using formula \eqref{eq:n3.7}.
The first one is  to use the local chart of $U_j$, i.e., $U_j$ will play 
the role of $U$ in \eqref{eq:n3.7}, and we get a function 
$u_{j,\epsilon}$. Similarly, the second way is to use the local chart 
of $U_k.$ Let $F$ be the change of coordinates on $U_j\cap U_k$  
from $U_j$ to $U_k.$
Denote by  $u_{j,\delta}^F$  the function given by 
Lemma \ref{lem_regularization} which corresponds to  
the regularization of $u_j$ using the local chart of  $U_k.$  Write
\begin{align*}
u_{j,\epsilon}-u_{k,\epsilon}=u_{j,\epsilon}-u_{j,\epsilon}^F 
+(u_j-u_k)_\epsilon\quad\text{on}\quad U_j\cap U_k,
\end{align*}
where  the  term $(u_j-u_k)_\epsilon$ is the regularization of 
$u_j-u_k$ by formula \eqref{eq:n3.7} using  the local chart of 
$U_k.$
Recall from  \eqref{eq:n3.16} that $u_j-u_k=f_j-f_k$ which is 
a smooth function.
This together with the previous equality and   
Lemma \ref{lem_regularization}, imply
\begin{align}\label{eq:n3.17}
\| (u_{j,\epsilon}-u_{k,\epsilon})-(f_j-f_k)\|_\infty
\lesssim M\epsilon^{\alpha}\quad\text{on}\quad U_j\cap U_k.
\end{align}
  
Fix  a  constant $c>0$ large  enough. For each $j\in J$ let $\eta_j$ 
be  a smooth function defined
in $U_j$ such that  $\eta_j=0$ on $V_j$ and that
$\eta_j=-c$  away from  a  compact subset of $U_j.$ 
We have that $\ddc \eta_j\geq  -c'\omega_0$  
for some constant $c'>0.$
 For each $\epsilon>0$ and $j\in J,$ consider the function 
 \begin{align}\label{eq:n3.18}
 v_j:= u_{j,\epsilon}-f_j+M\epsilon^\alpha \eta_j
 \quad \text{on}\quad U_j.
\end{align}
We identify  $J$ with $\{1,\ldots,l\}$   and set
\begin{align}\label{eq:n3.19}
 \phi_\epsilon:=M\epsilon^{\alpha-1} \maxep
\big(M^{-1}\epsilon^{1-\alpha} v_1, 
\ldots, M^{-1}\epsilon^{1-\alpha} v_l \big ).
\end{align}
Note that to define $\phi_\epsilon(x)$, $x\in X$,  
we remove $M^{-1}\epsilon^{1-\alpha}v_j$ from the last formula 
if $x\not\in U_j$. 

We first show that the function $\phi_\epsilon$ is
smooth on $X.$
For this purpose, we only  need to prove  the property in a 
neighborhood of an arbitrary fixed point of $X.$ Since  each $v_j$
is  well-defined and smooth on $U_j,$ 
using \eqref{eq:n3.19} and assertion  a) in  
Lemma \ref{lem_properties_regularized_max}, it is enough to prove 
the following claim. 

\medskip
\noindent
{\bf Claim 1.} For all $x\in U_j$  close enough to  $b U_j,$ we have
\begin{eqnarray*}
\lefteqn{\maxep \big(M^{-1}\epsilon^{1-\alpha} v_1, 
\ldots, M^{-1}\epsilon^{1-\alpha} v_l \big )(x)} \\
 & \qquad = &  \maxep
\big(M^{-1}\epsilon^{1-\alpha} v_1, \ldots, 
\widehat{M^{-1}\epsilon^{1-\alpha} v_j} , \ldots, M^{-1}
\epsilon^{1-\alpha} v_l \big )(x). 
\end{eqnarray*}

\medskip

Let  $k\in J$ such that $x\in V_k.$
We infer from  \eqref{eq:n3.18} and the  equality $\eta_k(x)=0$
 that
 $$  v_k(x) = u_{k,\epsilon}(x)-f_k(x) .$$
The  same  argument  using the  equality $\eta_j(x)=-c$ gives
$$  v_j(x) = u_{j,\epsilon}(x)-f_j(x)-cM\epsilon^\alpha. $$
Putting  the two last equalities  together  with \eqref{eq:n3.17},  
and using that $c>0$ is large  enough, we infer 
$$v_k(x)\geq  v_j(x)+2M\epsilon^\alpha.$$ 
 This, combined with  assertion c)  in  
Lemma \ref{lem_properties_regularized_max}, implies Claim 1.

\medskip
\noindent
{\bf Claim 2.} The function $\phi_\epsilon$ belongs to 
  $\PSH(X,(1+c'M\epsilon^\alpha) \omega_0)$. 

\medskip

It is enough to work in a small open set $W$ in $X$. By Claim 1, 
we can remove from the definition \eqref{eq:n3.19}  
of $\phi_\epsilon$ all functions 
$M^{-1}\epsilon^{1-\alpha} v_j$ if $W\not\subset U_j$.  
So we have $W\subset U_j$ for the indexes $j$ considered below.
Since $u_j$ is p.s.h., so is  $u_{j,\epsilon}.$
Therefore,  we deduce from  \eqref{eq:n3.15} and
\eqref{eq:n3.18} that 
\begin{align*}
\ddc v_j=\ddc u_{j,\epsilon}-\omega_0 +M\epsilon^\alpha 
\ddc \eta_j\geq -(1 +c'M\epsilon^\alpha)\omega_0.
\end{align*}
Choose a function $f$ on $W$ such that 
$\ddc f= M^{-1}\epsilon^{1-\alpha} (1+c'M\epsilon^\alpha)
\omega_0$. We deduce from \eqref{eq:n3.19} and the construction 
of $\max_\epsilon$ that 
$$\phi_\epsilon=M\epsilon^{\alpha-1} \maxep
\big(M^{-1}\epsilon^{1-\alpha} v_1+f, \ldots,
M^{-1}\epsilon^{1-\alpha} v_l +f\big)-M\epsilon^{\alpha-1}f.$$
Since $M^{-1}\epsilon^{1-\alpha} v_j+f$ is p.s.h. on $W$,  
applying assertion d) in
Lemma  \ref{lem_properties_regularized_max}, we obtain  
that $\phi_\epsilon$ belongs to 
  $\PSH(X,(1+c'M\epsilon^\alpha) \omega_0),$ thus proving Claim 2.

\medskip

We continue the proof of the theorem.  By \eqref{eq:n3.16} 
and \eqref{eq:n3.18},  we get on $V_j$
$$\|\phi-v_j  \|_\infty= \|(u_j-f_j)-(u_{j,\epsilon}-f_j
+M\epsilon^\alpha \eta_j)\|_\infty
\leq \|u_j-u_{j,\epsilon}\|_\infty+M\epsilon^\alpha \|\eta_j\|_\infty
\lesssim M\epsilon^{\alpha}.$$
 This and assertion b) in 
 Lemma \ref{lem_properties_regularized_max} prove the first 
 estimate in \eqref{eq:n3.12}.
For the second estimate,  we infer from assertion e) of 
Lemma  \ref{lem_properties_regularized_max} that 
\begin{eqnarray}\label{eq:n3.26}
\| \phi_\epsilon\|_{\Cc^m} &=&  M\epsilon^{\alpha-1}  
\left\|  \maxep ( M^{-1}\epsilon^{1-\alpha}v_1,\ldots, 
 M^{-1}\epsilon^{1-\alpha}v_l) \right\|_{\Cc^m} \nonumber\\
&\lesssim &   M\epsilon^\alpha+  \sup_{1\leq i\leq l} \| v_i\|_\infty 
+  M\epsilon^{\alpha-1}\sum_{r_{ij}}
  \epsilon^{1-\sum r_{ij}} \prod_{i,j}  \big(M^{-1}\epsilon^{1-\alpha}
  \| v_i\|_{\Cc^j}\big)^{r_{ij}},
\end{eqnarray}
 the sum being taken over all $r_{ij}>0$ with $1\leq i\leq l$ and 
 $j\geq 1$ such that $\sum jr_{ij}\leq  m.$  
On the  other hand,  by \eqref{eq:n3.9} and \eqref{eq:n3.18}, 
we have
\begin{align*}
\| v_i\|_{\Cc^j} =  \|u_{i,\epsilon}-f_i 
+M\epsilon^\alpha \eta_i\|_{\Cc^j}\lesssim  M\epsilon^{-j+\alpha}.
\end{align*}
Inserting these estimates  into  \eqref{eq:n3.26},  we obtain that
 $\phi_\epsilon$ satisfies  the second inequality in \eqref{eq:n3.12}. The theorem follows. Note that we can get similar estimates for every $m\in\R_+$. 
\hfill $\square$

\begin{remark}\label{rk_Demailly} \rm
We can prove in the same way the existence of  constants
$c> 0$ depending only  on $X,$ $\omega_0,$ and
$c_m>0$   depending only on 
$X,\omega_0,m\in\N^*$,  satisfying the following property.
Let $\phi$ be  an  $\omega_0$-p.s.h. function in $\LLip(X)$.
Then,  for each $0<\epsilon\leq 1/2,$ there  exists a smooth 
function $\phi_\epsilon$  such that
 \begin{enumerate}
 \item[a)] $\phi_\epsilon$ is $\omega_0$-p.s.h.; 
 \item[b)]  $\|\phi_\epsilon -\phi\|_\infty
 \leq -c (1+\|\phi\|_\LLip) \epsilon\log \epsilon;$
 \item[c)] $\| \phi_\epsilon\|_{\Cc^m(X)}
 \leq -c_m (1+\|\phi\|_\LLip) \epsilon^{-m+1} \log\epsilon $  
 for $m\in \N^*.$
 \end{enumerate}
\end{remark}

\subsection{Regularity of equilibrium weight}

In this subsection, we  study the equilibrium weight associated 
with a weighted compact subset $(K,\phi)$ of $X$. 
We start with the following tautological maximum  principle, 
and we refer the reader  to the beginning of the paper and
the Introduction for the notation used  below.

 \begin{proposition}  \label{prop_max_principle}
 Let  $ (K,\phi) $ be a regular weighted subset of $X$ and let 
 $P_K\phi$ be the associated equilibrium weight. Then 
 for every $\omega_0$-p.s.h. function $\psi$ on $X$, we have 
$$\sup_K(\psi-\phi)=\sup_K(\psi-P_K\phi)=\sup_X(\psi-P_K\phi).$$
In particular, for every section $s\in H^0(X,L^p)$ we have 
$$\|s\|_{L^\infty(K,p\phi)}=\|s\|_{L^\infty(K,pP_K\phi)}
=\| s\|_{L^\infty(X,pP_K\phi)}.$$
 \end{proposition}
 \proof 
 By Definition \ref{def_weighted_K}, we have $P_K\phi\leq \phi$ 
 on $K.$ Hence, 
$$\sup_K(\psi-\phi)\leq  \sup_K(\psi-P_K\phi)
\leq  \sup_X(\psi-P_K\phi).$$
 To prove the converse  inequality, observe that 
 $\psi- \sup_K(\psi-\phi) \leq \phi$ on $K.$ 
 This, combined with  Definition \ref{def_weighted_K}  and  the fact 
 that $\psi$ is $\omega_0$-p.s.h., 
 implies that $\psi-  \sup_K(\psi-\phi) \leq P_K\phi$ on $X.$  
 We deduce $\psi-P_K\phi\leq \sup_K (\psi-\phi)$ and then
 the first assertion in the proposition.
 
Next,  observe  that
$$\ddc  {1\over p} \log|s|={1\over p} [s=0] 
-\omega_0\geq -\omega_0,$$
where $[s=0]$ is the current of integration on the hypersurface 
$\{s=0\}$. So 
 ${1\over p} \log|s|$ is $\omega_0$-p.s.h.
 Applying the first assertion of the proposition to this function
 instead of $\psi$ gives the second assertion.
\endproof 

The following basic result  has been stated in  
\cite[Lemma 2.14]{BermanBoucksom}.

\begin{lemma} \label{lem_projection_P_K_is_Lipschitz}
Let $K$ be a non-pluripolar compact subset of $X$. Then 
the  projection $P_K$ is  non-decreasing, concave, and  
continuous along decreasing sequences of continuous weights 
$\phi$ on $K.$ It is also $1$-Lipschitz continuous, that is,
$$\sup_X|P_K\phi_1-P_K\phi_2|\leq \sup_K|\phi_1-\phi_2|$$
for all continuous weights $\phi_1$ and $\phi_2$ on  $K.$
\end{lemma}
\proof 
We only give  the proof of the inequality in the lemma and  leave  
the  verification of the  other  statements
to the  reader.  Since  $\phi_1\leq \phi_2 +  \sup_K|\phi_1-\phi_2|$ 
on $K,$ it follows  from  Definitions  \ref{def_weighted_K} and
\ref{def_K_regular} that
$$P_K\phi_1\leq  P_K\phi_2+  \sup_K|\phi_1-\phi_2|\quad
\text{on}\quad X.
$$
This and  the similar  estimate  which is obtained by interchanging
$\phi_1$ and $\phi_2,$ imply the desired inequality.
\endproof 

The following theorem is  the  main result of this  subsection. 
It gives us a class of compact sets $K$ satisfying regularity 
properties mentioned in the Introduction.

\begin{theorem}\label{th_reg_K}
Let $K$ be the closure of a non-empty open subset of $X$ with  
$\Cc^2$ boundary. Then $K$ is $(\Cc^\alpha,\Cc^{\alpha})$-regular  
for every $0<\alpha<1$. 
\end{theorem}

It is known that such a compact set is regular. To prove 
this property, it is enough to show that $P_K\phi$ is continuous 
when $\phi$ is 
H\"older continuous and then obtain the same property for 
continuous $\phi$ by approximation. Thus, the regularity of $K$ 
can be also obtained with the arguments given below. 

\medskip
 
\noindent{\bf Proof of Theorem \ref{th_reg_K} in the case $K=X$.} 
Let $\phi$ be a $\Cc^\alpha$ function on $X$ with bounded 
$\Cc^\alpha$-norm. 
We have to show that $\psi:=P_X\phi$ has bounded 
$\Cc^\alpha$-norm.
We will need to regularize $\psi$ using the method 
 introduced by Demailly in \cite{Demailly94}. 
Recall that for simplicity we use here the metric on $X$ induced 
by the K\"ahler form $\omega_0$.

Consider the  
exponential map associated with the Chern connection on the 
tangent bundle $TX$ of $X$. The {\it formal holomorphic part} 
of its Taylor expansion is denoted by
$$\exph: TX\to X \quad \text{with} \quad T_zX\ni\zeta\mapsto 
\exph_z(\zeta).$$ 
It is approximatively the part of the exponential map which 
is holomorphic in $\zeta$, see \cite{Demailly94} for details. 
Let $\chi:\ \R\to [0,\infty)$ be a smooth
function with support in $(-\infty,1]$ defined by
$$\chi(t):={\const\over  (1-t)^2}\exp{1\over t-1}\quad 
\text{for}\quad  t<1,
\quad \chi(t)=0\quad \text{for} \quad t\geq 1,$$
where the constant $\const$  is adjusted so that 
$\int_{|\zeta|\leq 1}\chi(|\zeta|^2)d\Leb(\zeta)=1$  
with respect  to the  Lebesgue measure $d\Leb(\zeta)$  on
$\C^n\simeq T_zX$. Fix a constant $\delta_0>0$ small enough.
Define 
\begin{align}\label{eq:n3.38}
\Psi(z,t):=\int_{\zeta\in T_zX}  \psi(\exph_z(t\zeta))
\chi(|\zeta|^2)d\Leb(\zeta)\qquad \text{for}\qquad (z,t)
\in X\times [0,\delta_0].
\end{align}

By  \cite{Demailly94}, there is  a constant $b>0$ such that  
the function $t\mapsto \Psi(z,t)+bt$ is increasing for $t$ in 
$[0,\delta_0]$. 
Observe also that $\Psi(z,0)=\psi(z)$. By definition, 
$\psi=P_X\phi$ is bounded by $\min\phi$ and $\max\phi$. 
The values of $\Psi(z,t)$ are averages of values of $\psi$. 
So $\Psi(z,t)$ is also bounded by the same constants $\min\phi$ 
and $\max\phi$.

Consider  for $c>0$  and  $\delta\in (0,\delta_0]$
 the  {\it Kiselman-Legendre  transform} 
 \begin{equation}\label{eq_Kiselman-Legendre_transform}
 \psi_{c,\delta}(z):=\inf_{t\in (0,\delta]}\Big( \Psi(z,t)
 +bt-b\delta - c\log{t\over\delta}\Big).
 \end{equation}
Since $t\leq\delta\leq\delta_0$, we see that $\psi_{c,\delta}$ 
is bounded below by $\min\phi-b\delta_0$ and taking $t=\delta$ 
we also see that $\psi_{c,\delta}$ is bounded above by $\max\phi$. 
 
Using a result by Kiselman, it is not difficult to show 
(see \cite{Demailly94}, see also \cite[Lemma 1.12]{BermanDemailly}) 
that $\psi_{c,\delta}$ is quasi-p.s.h. and 
$$\omega_0+\ddc\psi_{c,\delta}\geq-(ac+b\delta)\omega_0,$$
where $a>0$ is  a constant, see also \cite{Kiselman78, Kiselman94}. 
Therefore, we have 
$$\ddc{\psi_{c,\delta}\over 1+ ac+b\delta}+\omega_0
\geq 0\quad \textrm{for all} \quad c>0.$$
From now on, we take $c=\delta^\alpha$. 
We have seen that $\psi_{c,\delta}$ is bounded uniformly in 
$c, \delta$ for $c$ and $\delta$ as above. Hence,
\begin{equation}\label{eq_estimate_psi_c,delta}
\Big|{\psi_{c,\delta}\over 1+ac+b\delta}-\psi_{c,\delta}\Big|
\lesssim\delta^\alpha.
\end{equation}
For $t:=\delta$ we obtain from 
\eqref{eq_Kiselman-Legendre_transform} 
that 
$$\psi_{c,\delta}(z)\leq \Psi(z,\delta).$$
On the  other hand, we  deduce from
\eqref{eq:n3.38}  that the value of  $\Psi(z,\delta)$  is an average
of the values $\psi$ in the ball  $\B(z,A\delta)$ in $X$ for 
some constant $A$ depending only on $X$ and $\omega_0.$
Since $\psi\leq \phi$ and the $\Cc^\alpha$-norm of $\phi$ is 
bounded, we have 
$$  \Psi(z,\delta)  \leq  \phi(z)+O(\delta^\alpha).$$
This, coupled with \eqref{eq_estimate_psi_c,delta}, gives 
$${\psi_{c,\delta}\over 1+ac+b\delta} 
\leq \phi+ O(\delta^\alpha).$$
Since the left hand side is an $\omega_0$-p.s.h. function,  
the identity $\psi=P_K\phi$ implies
$${\psi_{c,\delta}\over 1+ ac +b\delta}
\leq\psi + O(\delta^\alpha).$$
Then, using that $c=\delta^\alpha$, we get
$$ \psi_{c,\delta} \leq \psi + O(\delta^\alpha).$$
This and  \eqref{eq_Kiselman-Legendre_transform} imply the 
existence of   $t_z\in (0,\delta]$  such that
 \begin{equation}\label{eq_rho_psi}
 \Psi(z,t_z)+bt_z\leq \psi(z) +c\log{t_z\over \delta}
 +O(\delta^\alpha).
 \end{equation}
Recall that the function $t\mapsto \Psi(z,t)+bt$  is  increasing 
and observe that its value  at  $t=0$ is equal to $\psi(z)$. So  
the last identity implies
 $$c\log {t_z\over\delta}+O(\delta^\alpha)\geq 0.$$
Therefore, since $c=\delta^\alpha$, we have
$\theta \delta\leq t_z\leq \delta$, where $0<\theta<1$ 
is a constant.  
By \eqref{eq_rho_psi} and using again that 
$t\mapsto \Psi(z,t)+bt$  is  increasing, we obtain
 \begin{align}\label{eq:n3.49}
\Psi(z,\theta \delta)-\psi(z)\leq O(\delta^\alpha).
\end{align}

Fix a point $z\in X$ and local coordinates in a neighborhood of $z$ 
so that the metric on $X$ coincides at $z$ with the standard metric 
given by the coordinates. The function $\psi$ is the difference
between a p.s.h. function $\psi'$ and a smooth function. 
In particular, $\Delta\psi-\Delta\psi'$ is smooth. Denote by $\mu$ 
the positive measure defined by $\Delta\psi'$. Consider the 
following quantity involving the mass of $\mu$ on the ball $\B(z,r)$
 $$\nu(r):={(n-1)! \over \pi^{n-1} r^{2n-2}} \|\mu\|_{\B(z,r)} 
 \quad \text{for} \quad  0<r\ll 1.$$ 
Note that if instead of $\mu$ we use the measure defined by 
$\Delta\psi$, then the last quantity is changed by a term $O(r^2)$.
So in the following computation, the use of $\Delta\psi'$
is equivalent to the one of $\Delta\psi$.  
The advantage of $\Delta\psi'$ is that by Lelong's theorem, 
the above function $\nu(r)$ is increasing. 

According to  \cite[(4.5)]{Demailly94} and using that $\chi$ is 
strictly positive on $[0,1)$, we have  the following  Lelong-Jensen 
type  inequality
\begin{eqnarray*}
 \Psi(z,t)-\psi(z) &=& \int_0^t{d\over d\tau} \Psi(z,\tau)
 d\tau \nonumber\\
&\geq& \int_0^t{d\tau\over \tau} \Big [\int_{\B(0,1)} 
\nu(\tau|\zeta|)\chi(|\zeta|^2)d\Leb(\zeta) -O(\tau^2)\Big] 
\nonumber\\
 &\geq&  \int_{t/2}^t{d\tau\over \tau} \Big [\int_{1/2<|\zeta|<3/4} \nu(\tau|\zeta|)\chi(|\zeta|^2)d\Leb(\zeta)\Big] - O(t^2) \nonumber\\
 &\gtrsim & \int_{t/2}^t \tau^{1-2n} \|\mu\|_{\B(z,\tau/2)} d\tau 
 - O(t^2) \nonumber\\
 &\gtrsim & t^{2-2n}\|\mu\|_{\B(z,t/4)} - O(t^2).
\end{eqnarray*} 
 Combining this and \eqref{eq:n3.49}, we obtain 
$$\|\mu\|_{\B(z,t)} \lesssim t^{2n-2+\alpha} \quad \text{for} 
\quad t\ll 1.$$
The estimate is uniform in $z\in X$. Applying
Lemma \ref{lem_Hoelder} below gives the result.
\hfill $\square$

\medskip

To complete the proof of Theorem  \ref{th_reg_K} for
$K=X$,  it remains  to prove the following elementary result, 
see also \cite{DinhNguyen}. For the reader's convenience, 
we give here a proof.

\begin{lemma}\label{lem_Hoelder}
Let $\phi$ be a subharmonic function in a neighborhood $U$ 
of $\overline{\B(0,1)}\subset \R^m$ and  
$0<\alpha< 1.$ Suppose there are  constants  $A>0$ and $t_0>0$  
such that  $\|\phi\|_\infty\leq A$, and 
for every $x\in \B(0,1)$ and $0<t\leq t_0$,  we have 
\begin{align}\label{eq:n3.55}
\|\Delta\phi\|_{\B(x,t)} \leq A t^{m-2+\alpha}.
\end{align}
Then $\phi$ is of class $\Cc^\alpha$ and its $\Cc^\alpha$-norm 
on $\B(0,1)$ is bounded by a constant depending only on
$U, A,t_0$ and $\alpha$. The result still holds for $\alpha=1$ 
if we replace $\Cc^\alpha$ by $\LLip$. 
\end{lemma}
\proof 
For simplicity, we only consider $0<\alpha< 1$ and $m\geq 3$.  
In this case, the Newton kernel $E(x)$ for $x\in \R^m$ is equal to 
a negative constant times $|x|^{2-m}$ and  
 $\Delta(E\ast \mu)=\mu$ for all measure $\mu$  with compact 
 support, see \cite[Theorem 3.3.2]{Hoermander}.
We can assume that $U=\B(0,1+4r_0)$ for some constant 
$r_0<t_0/4$ and that $\Delta\phi$ has finite mass in $U$. So 
\eqref{eq:n3.55} holds for $t\leq 4r_0$. 
Define $\mu:=\Delta \phi$ on $U$ 
and  $f:= E\ast \mu.$ The function $f-\phi$ is harmonic  on $U$. 
Therefore, we only 
need to show that  $f$ has bounded $\Cc^\alpha$-norm on 
$\B(0,1)$. 
  
Fix two points $x,y\in \B(0,1)$ and define $r:={1\over 2}|x-y|$. 
Since $\|\phi\|_\infty\leq A$, we only need to show that
$|f(x)-f(y)|\lesssim r^\alpha$ for $r\ll r_0$. 
Define 
$$D_1:=\B(x,r),\ \  D_2:=\B(y,r), \ \  
D_3:=\B(x,r_0)\setminus (D_1\cup D_2),  
D_4:=   \B(0,1+4r_0) \setminus \B(x,r_0)$$
and
$$I_k:=\int_{D_k}  \big ||x-z|^{2-m}- |y-z|^{2-m}\big|d\mu(z).$$
Observe that  $|f(x)-f(y)|\lesssim  I_1+I_2+I_3+I_4$. So it is enough
to bound $I_1, I_2, I_3,I_4$. 

Consider the integral $I_1$. The case of $I_2$ can be treated 
in the same way. 
Since $|z-x|\leq |y-z|$ for $z\in D_1$, we have 
\begin{align}\label{eq:n3.62}
I_1\leq 2 \int_{\B(x,r)}|x-z|^{2-m} d\mu(z). 
\end{align}
Recall that $\mu=\Delta\phi$ and it satisfies  \eqref{eq:n3.55}. 
Observe that $|x-z|^{2-m}$ can be bounded by a constant times
the following combination of the characteristic functions of balls
$$|x-z|^{2-m}\lesssim \sum_{k=0}^\infty (2^{-k}r)^{2-m}
\ind_{\B(x,2^{-k}r)}.$$
The integral in \eqref{eq:n3.62} is bounded by a constant times 
$$\sum_{k=0}^\infty  (2^{-k}r)^{2-m} \|\Delta\phi\|_{\B(x,2^{-k}r)}
\lesssim \sum_{k=0}^\infty \int_{2^{-k}r}^{2^{-k+1}r} \tau^{1-m} 
\|\Delta\phi\|_{\B(x,\tau)}d\tau=\int_0^{2r}\tau^{1-m} 
\|\Delta\phi\|_{\B(x,\tau)}d\tau.$$
We then deduce from \eqref{eq:n3.55} that $I_1\lesssim r^\alpha$.

Consider now the integral $I_3$. Observe that $|x-z|\approx  |y-z|$ 
when  $z\not\in D_1\cup D_2.$ Hence
\begin{align}\label{eq:n3.69}
\big| |x-z|^{2-m}  - |y-z|^{2-m}\big| \lesssim  r|x-z|^{1-m}
\end{align}
and 
$$I_3\lesssim r\int_{\B(x,r_0)\setminus \B(x,r)}
|x-z|^{1-m}d\mu(z).$$
We need to bound the last integral by $O(r^{\alpha-1})$
and we can assume that $x=0.$  
Observe that we have on the domain $r<|z|<r_0$, 
$${1\over |z|^{m-1}} \lesssim \sum_{k=-\log_2 r_0}^{-\log_2 r} 
(2^{-k})^{1-m} \ind_{\B(0,2^{-k})}.$$
Hence, we obtain the following inequalities which imply the desired 
estimate for $I_3$
$$\int_{r<|z|<r_0} {d\mu(z)\over |z|^{m-1}} 
\lesssim\sum_{k=-\log_2r_0}^{-\log_2r} (2^{-k})^{1-m}
\|\mu\|_{\B(0,2^{-k})} \lesssim \sum_{k=-\log_2 r_0}^{-\log_2 r} 
(2^{-k})^{\alpha-1}.$$

Finally, for the integral $I_4$ with $z\in D_4$, observe that 
\eqref{eq:n3.69} implies 
$$\big| |x-z|^{2-m}  - |y-z|^{2-m}\big| \lesssim  r.$$
The estimate $I_4\lesssim r$ follows immediately. 
This completes the proof of the lemma.
\endproof

We continue the proof of  Theorem \ref{th_reg_K}.  
We need the following lemma.
For $r>0$ and $w\in\C$, denote by $\D(w,r)$  the  disc  of center 
$w$ and radius $r$ in $\C$. 

\begin{lemma}\label{L:Poisson}
Let $\alpha>0$ be a constant. Let
$u$ be  a quasi-subharmonic function on a neighborhood of 
$\overline{\D(-1,3)}$  such that 
$\Delta u\geq -1$,  $u\leq  1$ on $\overline{\D(-1,3)}$ and 
$u(z)\leq  |z|^\alpha$ for all $z\in \overline{\D(1,1)}$. Then 
there is a constant $c>0$ depending only on $\alpha$ such that for 
all $t\in [-1/2,0]$ we have
$u(t)\leq  c|t|^{\min(1,\alpha)}$ if $\alpha\not=1$ 
and $u(t)\leq -c|t|\log|t|$ if $\alpha=1$. 
\end{lemma} 
 \proof 
 Replacing $\alpha$ by $\min(2,\alpha)$ allows us to assume that 
 $\alpha\leq 2$. 
Observe that the function $|z|^2$ is smooth and its Laplacian is 
equal to 2. So
replacing $u(z)$ by ${1\over 20}\big[{u(z)+|z|^2}\big]$ allows 
us to assume, from now on,  that $u$ is subharmonic.
 Let $\Omega$ denote the domain 
 $\D(-1,3)\setminus \overline{\D(1,1)}$.
Let $\Phi:\Omega\to\D(0,1)$ be a bi-holomorphic map which sends
$-4,0$ and $[-4,0]$ to $-1,1$ and $[-1,1]$, respectively. 
Since $b\Omega\setminus\{2\}$ is smooth analytic real, by 
Schwarz reflexion, $\Phi$ can be extended to a holomorphic map 
in a neighborhood of this curve and $\Phi'$ does not vanish there. 

Define $z'=\Phi(z)$ and $v(z'):=u\circ\Phi^{-1}(z')=u(z)$. 
We deduce from $u(z)\leq |z|^\alpha$ that 
$v(z')\lesssim |z'-1|^\alpha$ for $z'\in b\D(0,1)$.
Let $t$ be as in the statement of the lemma and define 
$t':=\Phi(t)$ and $s:=1-t'$. We have $s\in [0,2]$ and 
$s\lesssim  |t| \lesssim s$. We only have to show that 
$v(t')\lesssim s^{\min(1,\alpha)}$ if $\alpha\not=1$ and 
$v(t')\lesssim-s\log s$ if $\alpha=1$. Since $v$ is subharmonic,
it satisfies the following inequality involving the Poisson integral 
on the unit circle
$$v(t')\lesssim \int_{-\pi}^\pi {1-|t'|^2\over |e^{i\theta}-t'|^2}
v(e^{i\theta})d\theta.$$
Observe that $1-|t'|^2\lesssim s$ and 
$|e^{i\theta}-t'|^2\gtrsim s^2+\theta^2$. The last inequality is clear 
for $\theta<4s$ because $|e^{i\theta}-t'|\gtrsim s$ as $t'$ cannot 
be too close to $-1$, and it is also clear when $\theta\geq 4s$. 
We then deduce from the estimate of $v$ on the unit circle that
$$v(t')\lesssim \int_{-\pi}^\pi {s|\theta|^\alpha\over s^2
+\theta^2}d\theta=s^{\alpha}\int_{-\pi/s}^{\pi/s}
{|\theta'|^\alpha\over 1+\theta'^2}d\theta'\leq s^{\alpha}
\int_{-\infty}^{\infty} {|\theta'|^\alpha\over 1+\theta'^2}d\theta'.$$

When $\alpha<1$, the last integral is finite and the lemma follows. 
Using the integral before the last one, we also see that if
$\alpha=1$ then $v(t')\lesssim -s\log s$ which also implies the 
lemma in this case.
Consider now the case $\alpha>1$. We deduce from the above
inequality that
$$v(z)\lesssim  s\int_{-\pi}^\pi |\theta|^{\alpha-2} 
d\theta\lesssim s.$$
This completes the proof of the lemma. 
\endproof

\noindent
{\bf  Proof of Theorem \ref{th_reg_K} in the case $K\not=X$.} 
Consider a weight $\phi$ of bounded $\Cc^\alpha$-norm on $K$ 
with $0<\alpha<1$. 
Adding to $\phi$ a constant allows us to assume that $\phi\geq 0$.
Dividing $\phi$ and $\omega_0$ by a constant allows us to assume 
that $\|\phi\|_{\Cc^\alpha}\leq 1/100$. 
We have to show that $P_K\phi$ is of class $\Cc^{\alpha}$.

Fix a large constant $A\gg \|\phi\|_{\Cc^\alpha}$ and define
$$\widetilde\phi(x):=\min_{y\in K} \big[\phi(y)
+A\dist (x,y)^{\alpha}\big] \quad \text{for}\quad x\in X.$$
Since $\phi$ is $\Cc^\alpha$ and $A$ is large, $\widetilde\phi$ is 
an extension of $\phi$ to $X$, i.e., $\widetilde\phi=\phi$ on $K$. 
Moreover, if the above minimum is achieved at a point $y_0\in K$, 
by definition of $\widetilde\phi$,  we have for $x'\in X$
$$\widetilde\phi(x')-\widetilde\phi(x)\leq \big(\phi(y_0)
+A\dist(y_0,x')^{\alpha}\big) -\big(\phi(y_0)
+A\dist (y_0,x)^{\alpha}\big) \leq A\dist(x,x')^{\alpha}.$$
Therefore, the function $\widetilde\phi$ is $\Cc^{\alpha}$. 

The idea is to reduce the problem to the case $K=X$ which was 
already treated above. We only need to show that 
$P_K\phi\leq \widetilde\phi$ because this inequality implies that
$P_K\phi=P_X\widetilde\phi$. 
Moreover, since $P_K\phi$ is bounded and $A$ is large enough, 
we only need to check that $P_K\phi(x)\leq \widetilde\phi(x)$ for 
$x$ outside $K$ and close enough to $K$.  

Fix a finite atlas with  local holomorphic coordinates (that we always 
denote by $z=(z_1,\ldots,z_n)$) on open subsets $U_i$ of $X$ 
satisfying the following properties
\begin{enumerate}
\item Each open set $U_i$ corresponds to a ball $\B(a_i,10)$ of 
radius 10 centered at some point $a_i$ in $\C^n$;
\item If $V_i\subset U_i$ denotes the open set corresponding to 
$\B(a_i,1)$, then these $V_i$ cover $X$; 
\item $\phi$ restricted to $K\cap U_i$ is identified to a  function on 
a subset of $\B(a_i,10)$; we still
denote this function by $\phi$; it satisfies 
$\|\phi\|_{\Cc^\alpha}\leq 1/100$; 
for simplicity, $K\cap U_i$ will be also written as $K\cap \B(a_i,10)$;
\item $P_K\phi$ restricted to $U_i$ is identified  to a quasi-p.s.h. 
function on $\B(a_i,10)$ that we still denote by $P_K\phi$; 
it satisfies  $P_K\phi\leq \phi$ on $K\cap \B(a_i,10)$ and  
$\ddc P_K\phi\geq -\omega_0\geq -{1\over 2}\ddc\|z\|^2$ on 
$\B(a_i,10)$;
\item For any point $y$ in $bK\cap \B(a_i,2)$, $K$ contains a ball 
$B$ of radius 2 such that $y\in bB$ and $bB$ is tangent to $bK$ 
at $y$.  This can be done because $K$ has $\Cc^2$ boundary.
\end{enumerate}
This choice of atlas does not depend on $A$. So we can increase 
the value of $A$ when necessary.

Now,  $x$ belongs to some $V_i$. In what follows, we drop 
the index $i$ for simplicity, e.g. we will write $a$ instead of $a_i$.
Recall that the point $x$ is assumed to be outside and near the set 
$K$. Let $y_0$ be as above and denote by $x_0$ the projection 
of $x$ to the boundary of $K$, i.e., $|x-x_0|=\inf_{y\in K} |x-y|$. 
Here, we use the standard metric on $\C^n$. This point $x_0$ is 
unique because $K$ has $\Cc^2$ boundary and $x$ is close to $K$. 
Define $r:=|x-x_0|$ which is a small number.

\medskip
\noindent
{\bf Claim.} We have $|x_0-y_0|\lesssim r$ and hence 
$y_0\in \B(a,2)$ and $\widetilde\phi(x)
\geq \phi(x_0)+A'r^{\alpha}$, where $A'>0$ is a big constant 
(if we take $A\to\infty$ then $A'\to\infty$). 

\medskip

Indeed, if the first inequality were wrong, we would have 
$|x-x_0|\ll |x-y_0|\approx |x_0-y_0|$ and
by definition of $\widetilde\phi(x)$ and $y_0$
$$\widetilde\phi(x)=\phi(y_0)+A\dist(x,y_0)^{\alpha}
\leq \phi(x_0)+ A\dist(x,x_0)^{\alpha}.$$ 
Note that the distance on $U\subset X$ is comparable with 
the Euclidean distance with respect to the coordinates $z$. 
This comparison is independent of $A$. 
So the inequality implies
$$\phi(x_0)-\phi(y_0)\gg |x_0-y_0|^{\alpha}$$
which is a contradiction because $\phi$ is $\Cc^\alpha$. 

We also obtain the second inequality in the claim using the definition 
of $\widetilde\phi, y_0,x_0, r$ and the first inequality  
$$\widetilde\phi(x)-\phi(x_0)=\phi(y_0)-\phi(x_0)
+A\dist(x,y_0)^{\alpha}\gg r^\alpha,$$
since $A$ is large, $\phi$ is $\Cc^\alpha$, and 
$|x-y_0|\geq |x-x_0|=r$.

\medskip

By the claim, it is enough to show that
$P_K\phi(x)\leq \phi(x_0)+A'r^{\alpha}$. 
Using a unitary change of coordinates, we can assume that $x_0$ 
and $x$ 
are the points of coordinates $(0,0,\ldots,0)$ and $(-r,0,\ldots,0)$, 
respectively. This change of coordinates does not change the
metric on $\C^n$, so it does not change the norms of functions.
We use the coordinate $z_1$ in the complex line
$\Lambda:=\{z_2=\cdots=z_n=0\}$ and denote by $\D(w,r)$ 
the disc of center $w$ and radius $r$ in $\Lambda$. 

We will apply Lemma \ref{L:Poisson} to a suitable function $u$. 
Recall that $\|\phi\|_{\Cc^\alpha}\leq 1/100$, $K$ has $\Cc^2$
boundary, $x_0$ is the projection of $x$ to $K$  and $r$ is small 
enough. By the choice of the coordinates $(z_1,\ldots,z_n)$, 
the intersection $K\cap \Lambda$ contains $\overline\D(1,1)$, 
see property (5) above.
Denote by $u$ the restriction to $\Lambda$ of the function 
$P_K\phi-\phi(x_0)$. 
We deduce from the definition of $P_K\phi$ and the above 
properties of the coordinates $z$ that $u$ satisfies the hypotheses 
of Lemma \ref{L:Poisson}. Therefore,  
$u(x)\lesssim r^{\alpha}$ and hence
$P_K\phi(x)-\phi(x_0)\lesssim r^\alpha$. 
This completes the proof of the theorem.
\hfill $\square$

\medskip

Note that the idea of the proof still works if instead of the ball $B$
in the above point (5) we only have a solid right circular cone of
vertex $y$ and of a given size such that its axis is orthogonal at
$y$ to the boundary of $K$. This allows us to consider the situation 
where $K$ is the closure of an open set whose boundary is not 
$\Cc^2$. We then need a version of Lemma \ref{L:Poisson} for 
an angle at 0 instead of $\D(1,1)$. This angle is equal to the 
aperture of the above circular cone. If $\theta\pi$ denotes this 
angle, then $K$ is $(\Cc^\alpha,\Cc^{\theta\alpha})$-regular for
$0<\alpha<1$. In the case of $\Cc^1$-boundary for example, 
we can choose $\theta$ as any constant strictly smaller than 1. 
As mentioned in the Introduction, we don't try to develop the paper
in this direction. 
 We thank Ahmed Zeriahi for notifying us the reference 
 \cite{PawluckiPlesniak} where Pawlucki and Plesniak considered 
 a class of compact sets which may be 
 $(\Cc^\alpha,\Cc^{\alpha'})$-regular.

\subsection{Asymptotic behavior of Bergman functions} 
\label{section_Bergman}

Recall that $(L,h_0)$ is a holomorphic Hermitian line bundle  
on a projective manifold $X$ whose
first Chern form is $\omega_0$. The probability measure $\mu^0$ is  associated with the volume form $\omega_0^n$  as in the beginning of the paper. 
We will work later with Hermitian metrics which are not necessarily 
smooth nor positively curved. It is crucial to understand
the asymptotic behavior of the Bergman kernel associated with
$L^{p}$ and the new metrics when $p$ tends to infinity. 

As mentioned above, our strategy is to approximate the considered 
metrics by smooth positively 
curved ones. So we need to control the dependence of the 
Bergman kernels in terms of the positivity of the curvature. 
The solution to this problem will be presented below.
 We refer to \cite{MaMarinescu07} for basic properties of 
 Bergman kernel.

Consider a metric $h=e^{-2\phi}h_0$ on $L,$ where  $\phi$ 
is a  continuous weight on a compact subset $K$ of  $X$. 
Recall that $H^0(X,L^p)$ denotes the space of holomorphic sections 
of $L^p$. Since $L$ is ample,  by Kodaira-Serre vanishing and 
Riemann-Roch-Hirzebruch theorems 
(see \cite[Thm 1.5.6 and 1.4.6]{MaMarinescu07}) we have
\begin{equation}\label{e:Demailly}
N_p:=\dim H^0(X,L^p) =  {p^n\over n!} \|\omega_0^n\| +O(p^{n-1}).
\end{equation}
Let  $\mu$ be a probability measure with support in $K$.  
Consider the natural $L^\infty$ and $L^2$ semi-norms on 
$H^0(X,L^p)$ induced by the metric $h$ on $L$ and the measure 
$\mu$, which are defined for $s\in H^0(X,L^p)$ by
\begin{align}\label{eq:n1.15}
\|s\|_{L^\infty(K,p\phi)}:=\sup_K|s|_{p\phi} \qquad \text{and} 
\qquad 
\|s\|_{L^2(\mu,p\phi)}^2:=\int_X |s|_{p\phi}^2 d\mu.
\end{align}

We will only use measures $\mu$ such that the above semi-norms 
are  norms, i.e., there is no section 
$s\in H^0(X,L^p)\setminus \{0\}$ which vanishes on $K$ or 
on  the support of $\mu$. 
The first semi-norm is a norm when $K$ is not contained in 
a hypersurface of $X$. The second one is a norm
when $\mu$ is the normalized Monge-Amp\`ere measure 
with continuous potential because such a measure has no mass 
on hypersurfaces of $X$. This is also the case for any
Fekete measure of order $p$ as can be easily deduced from 
Definition \ref{def_Fekete_L}.

From now on, assume  that the above semi-norms are  norms and 
for the rest of this section, consider $K=X$. 
Let  $\{s_1,\ldots, s_{N_p}\}$ be an orthonormal basis of
$H^0(X,L^{p})$ with respect to the above $L^2$-norm.

\begin{definition} \rm \label{defi_Bergman_distortion} 
We call   {\it Bergman function of $L^{p}$}, 
associated with $(\mu,\phi)$, the function $\rho_p(\mu,\phi)$  
on $X$ given by
$$\rho_p(\mu,\phi)(x):=\sup\Big\{|s(x)|_{p\phi}^2: \ \ 
s\in H^0(X,L^p), \|s\|_{L^2(\mu,p\phi)}=1\Big\} 
=\sum_{j=1}^{N_p}|s_{j}(x)|_{p\phi}^2$$
and we define {\it the Bergman measure} associated with  
$(\mu,\phi)$ by
$$\Bc_p(\mu,\phi):=N_p^{-1}\rho_p(\mu,\phi)\mu.$$
\end{definition}
Note that it is not difficult to obtain the identity in the definition of 
$\rho_p(\mu,\phi)$ and  check that $\Bc_p(\mu,\phi)$ is a 
probability measure.  For the above definition, we only need that 
$\phi$ is defined on the support of $\mu$ or a compact set 
containing this support.

\smallskip

In the rest of this subsection, we  assume that the weight $\phi$
is a function of class $\Cc^3$ on $X$ and the first Chern form
 $\omega:=\ddc\phi+\omega_0$ satisfies 
 \begin{equation}\label{eq_strict_positivity}
\omega\geq \zeta\omega_0 \quad \text{for some constant}\quad 
\zeta>0.
\end{equation}
Note that this inequality implies that $\zeta\leq 1$ because 
$\omega$ and $\omega_0$ are cohomologous.
Here is the main result in this section which gives us an estimate of 
the Bergman function in terms of $\phi,\omega,p$ and $\zeta$.
We refer to the beginning of the paper for the notation.

\begin{theorem}\label{tue4} 
There exists  a constant $c>0$,  depending only on $X,L$ and the 
$\Cc^3$-norm of the Hermitian metric $h_0$ of $L$, with the 
following property.
For every $p>1$ and  every  weight $\phi$ 
of class $\Cc^{3}$ such that  \eqref{eq_strict_positivity} holds for 
some $\zeta$ with $\zeta\geq \|\phi\|_3^{2/3}(\log p)p^{-1/3}$, 
we have 
$$ \Big\|{\rho_p(\mu^0,\phi)(x)\over N_p}
- {\omega(x)^n\over\omega_0(x)^n}\Big\|_{L^1(\mu^0)} 
\leqslant  c\,  \|\phi\|_3 \zeta^{-3/2} (\log p)^{3/2}p^{-1/2}$$
with $\mu^0:=\|\omega_0^n\|^{-1}\omega_0^n$ 
the normalized Lebesgue measure on $X$, and
$$ \int_X\big | \Bc_p(\mu^0,\phi)(x)  - \mu_\eq(X,\phi)(x)\big |
   \leqslant c\,  \|\phi\|_3 \zeta^{-3/2} (\log p)^{3/2}p^{-1/2}.$$
\end{theorem}
\proof  
By hypotheses, $\phi$ is  $\omega_0$-p.s.h. Hence,
we have $\phi=P_X\phi$  and 
$\mu_\eq(X,\phi)=\NMA(\phi)=\|\omega_0^n\|^{-1}\omega(x)^n.$
Therefore, the second assertion is a direct consequence of 
the first one and Definition \ref{defi_Bergman_distortion}.

Consider now the first assertion. We use some  ideas  
from Berndtsson \cite[Sect. 2]{Berndtsson} and the recent joint 
work of Coman, Marinescu and the second author 
\cite{ComanMaMarinescu}, see also \cite{DinhMaMarinescu}.
Consider a point $x\in X$. Choose a local system of coordinates 
$z=(z_1,\ldots,z_n)$ centered at $x$ and a constant $c>0$ 
such that 
\begin{enumerate}
\item Some neighborhood of $x$ can be identified to the unit 
polydisc $\D^n$ in $\C^n$;
\item $\big\|\omega_0(z) - {\sqrt{-1}\over \pi }\sum_{j=1}^n dz_j 
\wedge d\bar z_j\big\| \leq c|z|$ for $z\in\D^n$;
\item   $\big|\phi(z) - q(z)-\sum_{j=1}^n (\lambda_j-1) |z_j|^2
\big|\leq c\|\phi\|_3 |z|^3$ for $z\in\D^n$, 
where $\lambda_i$ are real numbers and $q(z)$ is a harmonic 
polynomial in $z,\overline z$ of degree $\leq 2$.
\end{enumerate}

Observe that after choosing $z$ satisfying (1)-(2), we can take 
$q(z)$ as the harmonic part in the Taylor expansion of order 2 of 
$\phi$ at $x\equiv 0$; then, using a {\it unitary} change of 
coordinates allows us to assume that the non-harmonic part in
this Taylor expansion is given by a diagonal matrix.
So we have (1)-(3) and furthermore, the constant $c$ is controlled 
by the $\Cc^3$-norm of the metric $h_0$ on $L$.  
The numbers $\lambda_j$ and the coefficients of $q(z)$ can be 
controlled by the $\Cc^2$-norm of $\phi$. Note that if the metric 
$h_0$ of $L$ is $\Cc^4$, thanks to a standard property in K\"ahler 
geometry, we can replace $c|z|$ in (2) by $c|z|^2$. 

\medskip\noindent
{\bf Claim.} There is a holomorphic frame $\e$ of $L$ over $\D^n$ 
such that if $\phi_0:=-\log|\e|$ (see the beginning of the paper
for the notation), then
$$\Big|\phi_0(z)-\sum_{j=1}^n |z_j|^2\Big|\leq c|z|^3,$$ 
where $c>0$ is a constant depending only on $X,L$ and 
the $\Cc^3$-norm of $h_0$.

\medskip

We first prove the claim. Consider a frame $\widetilde \e$ of $L$ 
over $\D^n$. It can be chosen in a fixed finite family of local frames  
of $L$ over a finite covering of $X$. Define 
$\widetilde\phi_0:=-\log|\widetilde \e|$. We have by definition 
of curvature that $\omega_0=\ddc \widetilde\phi_0$. 
As above, thanks to (3), we can write 
$\widetilde\phi_0(z)=\widetilde q_0(z)
+\sum_{j=1}^n |z_j|^2+O(|z|^3)$, where $\widetilde q_0(z)$
is a harmonic polynomial of degree $\leq 2$. So we can write 
$\widetilde q_0(z)=\Re \widetilde Q_0(z)$, where 
$\widetilde Q_0(z)$ is a holomorphic polynomial of degree 
$\leq 2$ whose coefficients are controlled by the 
$\Cc^2$-norm of $h_0$. Define 
$\e=e^{\widetilde Q_0}\widetilde \e$. We have
$$|\e(z)|^2=|\widetilde \e(z)|^2 e^{2\widetilde q_0(z)}
=e^{2\widetilde q_0(z)-2\widetilde\phi_0(z)}.$$
The claim follows.

\medskip

Now, by (2) and (3), we have
$$ \omega(x)=\ddc\phi(x)+\omega_0(x)= {\sqrt{-1}\over \pi}  
\sum_{j=1}^n \lambda_j dz_j \wedge d\bar z_j.$$
 Hence, we get 
 \begin{align}\label{e:omega_vs_omega_0}
 \omega^n(x)=  \lambda_1\cdots\lambda_n\omega_0^n(x).
 \end{align}
Moreover, the inequality \eqref{eq_strict_positivity} at the point $x$ becomes
$$\lambda_j\geq \zeta \quad \text{for}\quad  1\leq j\leq n.$$
 Define 
\begin{align} \label{e:psi}
\varphi(z):=\sum_{j=1}^n \lambda_j |z_j|^2 \quad \text{and}\quad 
 \psi(z):=\phi(z) - q(z)-\varphi(z)+\phi_0(z).
 \end{align}
 
Consider  a normalized section $s\in H^0(X, L^p)$ with 
$\|s\|_{L^2(\mu^0,p\phi)}=1$. We are going to bound 
$|s(x)|_{p\phi}$ from above. Writing $s = f\e^{\otimes p},$ 
where $f$ is a holomorphic  function on  $\D^n$ and $\e$
is the frame given by the above claim. 
We apply the submean inequality for the p.s.h. function 
$|f(z)|^2e^{-2pq(z)}$ on the polydisc 
$\D_r^n:=\D_r\times\cdots\times \D_r$ ($n$ times) with radius 
$r:=(\log p)^{1/2}p^{-1/2}\zeta^{-1/2}$. Thanks to the special
form of $\varphi$, we obtain
\begin{align} \label{e:BM:s}
|s(x)|^2_{p\phi} =|f(0)|^2 e^{-2pq(0)}\leq {\int_{\D_r^n} 
|f|^2e^{-2pq-2p\varphi}   d\Leb\over \int_{\D_r^n} e^{-2p\varphi} 
d\Leb}\cdot
\end{align}
Note that the hypothesis on $\zeta$ and the fact that 
$\zeta\leq 1$ insure that $r\leq p\|\phi\|_3r^3\leq 1$. 
We will use this property in the computation below.

For the first integral in  \eqref{e:BM:s}, observe that by (2), 
the Lebesgue measure in $\D^n$ is equal to 
${1\over n!}\big({\pi\over 2}\big)^n\omega_0^n+O(|z|)$.
This, together with (3), \eqref{e:psi} and the above claim, gives
\begin{eqnarray*}
\int_{\D_r^n}   |f|^2e^{-2pq-2p\varphi}  d\Leb  &\leq &
\Big[{1\over n!} \Big({\pi\over 2}\Big)^n+O(r)\Big] \int_{\D_r^n}  
|f|^2e^{-2pq-2p\varphi} \omega_0^n\\
&\leq& \Big[{1\over n!} \Big({\pi\over 2}\Big)^n+O(r)\Big] 
\exp\big(2 p\max_{\D_r^n} \psi\big)\int_{\D_r^n}  
|f|^2e^{-2p(q+\varphi+\psi)}   \omega^n_0\\
&\leq &\Big[{1\over n!} \Big({\pi\over 2}\Big)^n+O(r)\Big] 
e^{O(p\|\phi\|_3r^3)}  \int_X  |s|^2_{p\phi}  \omega^n_0\\
&=& {1\over n!} \Big({\pi\over 2}\Big)^n \|\omega_0^n\|
+ O\big(\|\phi\|_3\zeta^{-3/2}(\log p)^{3/2} p^{-1/2}\big),
\end{eqnarray*}
because $\|s\|_{L^2(\mu^0,p\phi)}=1$ and 
$e^{O(p\|\phi\|_3r^3)}=1+O(p\|\phi\|_3 r^3)$. 

Define
$$E(t) :=\int_{\xi\in \D_t} e^{-2|\xi|^2} d\Leb(\xi) 
={\pi\over 2}(1-e^{-2t^2})\leq {\pi\over 2}\cdot$$
A direct computation shows that the second integral in 
\eqref{e:BM:s} is equal to
$$\int_{\D_r^n} e^{-2p\varphi}  d\Leb=\prod_{j=1}^n 
\int_{z_j\in \D_r} e^{-2p\lambda_j|z_j|^2} d\Leb(z_j)
=\prod_{j=1}^n {E(r\sqrt{p\lambda_j})\over p\lambda_j}\geq 
\Big({\pi\over 2}\Big)^n {(1-1/p^2)^n\over p^n
\lambda_1\ldots\lambda_n}$$  
since $r^2p\lambda_j\geq r^2p\zeta=\log p$. 

Combining the above estimates with \eqref{e:BM:s}, we obtain
$$ |s(x)|^2_{p\phi}\leq \Big[1
+O\big(\|\phi\|_3\zeta^{-3/2}(\log p)^{3/2} p^{-1/2}\big)\Big]
{1\over n!} p^n \lambda_1\ldots\lambda_n\|\omega_0^n\|.$$
By Definition \ref{defi_Bergman_distortion}, we get
$${ \rho_p(\mu^0,\phi)(x)\over p^n}\leq \Big[1
+O\big(\|\phi\|_3\zeta^{-3/2}(\log p)^{3/2} p^{-1/2}\big)\Big]
{1\over n!} \lambda_1\ldots\lambda_n\|\omega_0^n\|.$$
Then, using \eqref{e:Demailly} and \eqref{e:omega_vs_omega_0}, 
we  obtain
\begin{equation}\label{smp1.5}
{\rho_p(\mu^0,\phi)(x)\over N_p}\leq   
\big(1+c\|\phi\|_3 \zeta^{-3/2} (\log p)^{3/2}p^{-1/2}\big)
{\omega(x)^n\over\omega_0(x)^n} \quad \text{with} \quad c>0.
\end{equation}

Now,  define for simplicity 
$$\vartheta_1(x):={\rho_p(\mu^0,\phi)(x)\over N_p}, 
\quad \vartheta_2(x):={\omega(x)^n\over \omega_0(x)^n} \quad 
\text{and}\quad \epsilon:= c\|\phi\|_3\zeta^{-3/2}(\log p)^{3/2} 
p^{-1/2}.$$ 
So $\vartheta_1$ and $\vartheta_2$ are two positive functions 
of integral 1 with respect to the probability measure $\mu^0$. 
Inequality \eqref{smp1.5} says that 
$\vartheta_1\leq (1+\epsilon)\vartheta_2$. We need to check that 
$\|\vartheta_1-\vartheta_2\|_{L^1(\mu^0)}\lesssim \epsilon$. 
By triangle inequality, it is enough to check that 
$\|\vartheta_1-(1+\epsilon)\vartheta_2\|_{L^1(\mu^0)}
\lesssim \epsilon$. But since the function 
$\vartheta_1-(1+\epsilon)\vartheta_2$ is negative, it suffices 
to check that the integral of this function with respect to $\mu^0$ 
is larger than or equal to $-\epsilon$. 
A direct computation shows that this 
integral is in fact equal to $-\epsilon$. The proof of the theorem is 
now complete.
 \endproof

\section{Equidistribution of Fekete points}  
\label{section_equidistribution}

In this section, we will give the proofs of the main results stated 
in the Introduction. The estimates obtained in the previous section
allow us to use the strategy by Berman, Boucksom and Witt 
Nystr\"om. We refer to the beginning of the article for the notation.

\subsection{Energy, volumes and Bernstein-Markov property} 
\label{subsection_preliminaries}

Recall from \cite{BermanBoucksom} that the
{\it Monge-Amp\`ere energy functional} 
$\mE $, defined on bounded weights in 
$\PSH(X,\omega_0)$, is characterized by
$$\left. {d\over dt}\right|_{t=0} \mE ((1-t)\phi_1+t\phi_2)
=\int_X(\phi_2-\phi_1) \NMA(\phi_1).$$
So $\mE $ is only defined up to an additive constant, but the 
differences such as 
$\mE(\phi_1) -  \mE(\phi_2)$
are well-defined,  see also \eqref{eq_normalization}. 

Consider a non-pluripolar compact set $K\subset X$ and 
a continuous weight $\phi$ on $K$. Define the 
{\it energy at the equilibrium weight} of $(K,\phi)$ as
$$\mE_\eq(K,\phi):= \mE(P_K\phi).$$
This functional  is also  well-defined up to an additive constant. 
We will need the following property which was  
established in \cite[Th. B]{BermanBoucksom}.

\begin{theorem}\label{thm_differentiability}
The map  $\phi\mapsto \mE_\eq(K,\phi),$  
defined on the affine space of continuous weights on $K$,
is concave  and G\^ateaux differentiable, with directional derivatives 
given by integration against the equilibrium measure:
$$\left . {d\over dt}\right|_{t=0}  \mE_\eq (K,\phi+tv)
=\big\langle v ,\mu_\eq(K,\phi)  \big\rangle \quad
\text{for every continuous function } v \text{ on } K.$$ 
In particular, for all continuous weights $\phi_1$ and $\phi_2$ 
on $K$, we have
$$|\mE_\eq(K,\phi_1)-\mE_\eq(K,\phi_2)|
\leq \|\phi_1-\phi_2\|_\infty.$$
\end{theorem}

Note that the second assertion is obtained by taking the integral 
on $s\in [0,1]$ of the first identity applied to $\phi:=\phi_1+sv$ 
and $v:=\phi_2-\phi_1$. We use here the fact that 
$\mu_\eq(K,\phi)$ is a probability measure.

Let $\mu$ be a probability measure on $X$ and $\phi$ a continuous 
function on the support of $\mu$. The semi-norm
 $\|\cdot\|_{L^2(\mu,p\phi)}$ on $H^0(X,L^{p})$ is defined as
 in \eqref{eq:n1.15} and recall that we only consider measures 
 $\mu$ for which this semi-norm is a norm.  
Let $\mathcal B^2_p(\mu,\phi)$ denote the unit ball in  
$H^0(X,L^{p})$ with respect to this norm and 
$N_p:=\dim H^0(X,L^{p})$.
Recall from \cite {BermanBoucksom} the following 
$\mL_p$-functional
\begin{align}\label{eq:n4.5}
\mL_p(\mu,\phi):=  {1\over  2pN_p}  
\log\vol\mathcal B^2_p(\mu,\phi). 
\end{align}
Here, $\vol$ denotes the Lebesgue measure on the 
vector space $H^0(X,L^{p})$ which is  only defined up 
to a multiplicative constant.  
Note that  the  differences such as  
$\mL_p(\mu_1,\phi_1)- \mL_p(\mu_2,\phi_2)$ is well-defined 
and do  not depend on the choice of $\vol$ for any
 probability measures $\mu_1$ and $\mu_2$,  
 see also \eqref{eq_normalization}.  The functional 
 $\mL_p$ satisfies the following concavity property, 
 see  \cite[Proposition 2.4]{BBW}.

\begin{lemma} \label{lem_L_concave}
The  functional  $\phi\mapsto \mL_p(\mu, \phi)$ is concave on
 the space of all   continuous  weights on the support of $\mu$.
\end{lemma}

Recall from Definition \ref{defi_Bergman_distortion} that 
the Bergman measure $\Bc_p(\mu,\phi)$ is a probability measure. 
Note that when $\mu$ is the average of $N_p$ 
generic Dirac masses (more precisely, for points 
$x_1,\ldots,x_{N_p}$ such that 
the vector $\det (s_i(x_j))$ in the Introduction does not vanish), 
one can easily deduce 
from Definition \ref{defi_Bergman_distortion} that 
$\Bc_p(\mu,\phi)=\mu$, by considering sections vanishing on 
$\supp(\mu)$ except at a point. Such sections exist because 
$N_p=\dim H^0(X,L^p)$. This property holds in particular for 
Fekete measures of order $p$. 

The following relation between the functional $\mL_p(\mu,\cdot)$ 
and $\Bc_p(\mu,\cdot)$ has been established 
in \cite[Lemma 5.1]{BermanBoucksom}, see also 
\cite[Lemma 5.1]{Berndtsson2} and \cite[Lemma 2]{Donaldson}.

\begin{lemma}\label{lem_directional_derivative_L}
The directional derivatives of $\mL_p(\mu,\cdot)$ at 
a continuous weight $\phi$ on the support of $\mu$ are given 
by the integration against
the Bergman measure $\Bc_p(\mu,\phi),$ that is,
$$\left.{d\over  dt}  \mL_p(\mu, \phi+tv)\right |_{t=0}
= \langle v,\Bc_p(\mu,\phi)\rangle ,\quad  \text{with } v,\phi 
\text{ continuous on the support of } \mu.$$
In particular, for all continuous functions $\phi_1$ and $\phi_2$ on 
the support of $\mu$, we have
$$|\mL_p(\mu,\phi_1)-\mL_p(\mu,\phi_2)|
\leq \|\phi_1-\phi_2\|_\infty.$$
\end{lemma}

Note that as in Theorem \ref{thm_differentiability}, 
the second assertion of the last lemma is a direct consequence 
of the first one. 

Consider the norm
 $\|\cdot\|_{L^\infty(K,p\phi)}$ on $H^0(X,L^{p})$ defined  
 in \eqref{eq:n1.15}. 
Let $\mathcal B^\infty_p(K,\phi)$ denote the unit ball in 
$H^0(X,L^{p})$ with respect to this norm. Define 
\begin{align}\label{eq:n4.8}
\mL_p(K,\phi):=  {1\over  2pN_p}  
\log\vol\mathcal B^\infty_p(K,\phi).
\end{align}
We have the following elementary lemma. 

\begin{lemma}\label{lemma_L_K_mu}
If $\mu$ is a probability measure with $\supp( \mu)\subset K,$ 
then 
$$\mL_p(K,\phi)\leq   \mL_p(\mu,\phi).$$
\end{lemma}
\proof  
Since $\mu$ is a probability measure, 
we see that 
\begin{align}\label{eq:n4.11}
\| s\|_{L^2(\mu,p\phi)}\leq  \| s\|_{L^\infty(K,p\phi)},
\qquad   s\in H^0(X, L^{p}).
\end{align}
The lemma follows. 
\endproof 

We have the following property that we will only use in the case of 
$\omega_0$-p.s.h. weights.

\begin{lemma}\label{lemma_L-K}
 Let $\mu$ be a  probability measure and $K\subset X$ a compact 
 set with $\supp( \mu)\subset K.$
Assume the following strong Bernstein-Markov inequality: 
there  exists a constant $B>0$ such that
$$\sup_K \rho_p(\mu,\phi)\leq Bp^B \quad \text{for} \quad p>1.$$
Then there exists $c>0$ depending only on $B$  such that for 
$p> 1,$   we have
$$0\leq  \mL_p(\mu,\phi)-\mL_p(K,\phi)\leq  cp^{-1}\log p.$$
\end{lemma}

\proof   For all $p>1$ and section
 $s\in  H^0(X,L^{p}),$ by \eqref{eq:n4.11} and  
 Definition \ref{defi_Bergman_distortion},  we have  
\begin{align}\label{eq:n4.13}
\|  s\|_{L^2(\mu,p\phi)}  \leq   \|  s\|_{L^\infty(K,p\phi)}
\leq  e^{pc_p}  \|  s\|_{L^2(\mu,p\phi)} , 
\end{align}
where
\begin{align}\label{eq:n4.14}
c_p:= {1\over 2p}  \log\sup_K \rho_p(\mu,\phi).
\end{align}
Since the volume form $\vol$ is homogeneous of degree 
$2N_p=\dim_\R H^0(X,L^{p})$, it follows from \eqref{eq:n4.13} that
\begin{align*}
0\leq   \log {\vol\mathcal B^2_p(\mu,\phi)\over  
\vol\mathcal B^\infty_p(K,\phi)}\leq  2pN_pc_p .
\end{align*}
Hence, by  definition of the $\mL$-functionals in \eqref{eq:n4.5} 
and  \eqref{eq:n4.8}, we have
$$0\leq   \mL_p(\mu,\phi)-\mL_p(K,\phi)
={1\over 2pN_p}\log {\vol\mathcal B^2_p(\mu,\phi)\over  
\vol\mathcal B^\infty_p(K,\phi)}\leq c_p. $$
This, \eqref{eq:n4.14} and the assumed strong Bernstein-Markov 
inequality imply the lemma.
\endproof

The following result gives us a class of compact sets $K$ 
satisfying the strong Bernstein-Markov inequality  
stated in Lemma \ref{lemma_L-K} for $(X,P_K\phi)$ instead of 
$(K,\phi)$, see also \cite[section 1.2]{BBW}. 
We refer to the beginning of the article for the definition 
of $\mu^0$. 

 \begin{theorem} \label{th_BM_open}
Let $A>0$ and $\alpha,\alpha'>0$ be constants. Let $K\subset X$ 
be a $(\Cc^\alpha,\Cc^{\alpha'})$-regular compact set. 
Let $\phi$ be a function on $K$ such that 
$\|\phi\|_{\Cc^\alpha}\leq A$. Then there is a constant $B>0$ 
depending only on $X,L,h_0,K,A, \alpha$ and $\alpha'$ such that 
$$\sup_X\rho_p(\mu^0,P_K\phi)\leq Bp^B 
\quad \text{for}\quad p>1.$$
In particular, the statement holds when $K$ is the closure of 
an open set in $X$ with $\Cc^2$ boundary, $0<\alpha'<1$, 
$\alpha\geq \alpha'$ and $A>0$.
\end{theorem}
\proof  
The second assertion is a consequence of the first one and 
Theorem \ref{th_reg_K}. We prove now the first assertion.

It is enough to consider the case where $0<\alpha'< 1$. 
Since $K$ is $(\Cc^\alpha,\Cc^{\alpha'})$-regular,  the function  
$\psi:=P_K\phi$ has bounded $\Cc^{\alpha'}$-norm on $X$. 
Consequently, we only need to prove that 
 \begin{equation}\label{eq_distortion_Holder_weights}
\sup_X \rho_p(\mu^0,\psi)\lesssim p^{2n/\alpha'} \quad 
\text{for} \quad  p>1.
  \end{equation}
For this purpose,  fix a point $x\in X$ and a section 
$s\in H^0(X,L^p)$ such that $\|s\|_{L^2(\mu^0,p\psi)}=1$. 
By Definition \ref{defi_Bergman_distortion},  
it is enough to prove the estimate
\begin{equation}\label{eq_extremal_Bergman_function_estimate}
|s(x)|^2_{p\psi}\lesssim p^{2n/\alpha'}
\end{equation}
uniformly in $x$ and $s$. 

Choose local coordinates $z$ near $x$ such that $z(x)=0$ and for 
simplicity we still write $\psi(z)$ for the restriction of $\psi$ to 
a neighborhood of $x$. 
Fix also  a local holomorphic frame $\e$ of $L$ over 
a neighborhood of $x$ such that 
$|\e(0)|_{\psi}=e^{-\psi(0)}$. We can write 
$s(z)=f(z)\e^{\otimes p}(z)$, where $f(z)$ is a holomorphic function such that 
$|f(0)|e^{-p\psi(0)}=|s(0)|_{p\psi}$. 
So we need to check that $|f(0)|^2e^{-2p\psi(0)}
\lesssim p^{2n/\alpha'}$. Write  $\psi_\e(z):=-\log |\e(z)|_{\psi}$. 
This function differs from $\psi(z)$ by a pluriharmonic function. 
Therefore, it is also of class $\Cc^{\alpha'}$ and by definition 
we have $\psi_\e(0)=\psi(0)$. It follows that 
$|\psi_\e(z)-\psi(0)|\lesssim |z|^{\alpha'}$, and hence
\begin{eqnarray} \label{eq:n4.19}
\qquad p^{2n/\alpha'} \ = \ p^{2n/\alpha'} 
\| s \|^2_{L^2(\mu^0,p\psi)}
& \gtrsim &  p^{2n/\alpha'} \int_{|z|<{p^{-1/\alpha'}}} 
|f(z)|^2e^{-2p\psi_\e(z)} d\Leb(z) \\
& \gtrsim & p^{2n/\alpha'} \int_{|z|<{p^{-1/\alpha'}}}  
|f(z)|^2e^{-2p\psi(0)} e^{-cp|z|^{\alpha'}}  d\Leb(z)\nonumber
\end{eqnarray}
for some constant $c>0$.
 
Using the submean property for $|f(z)|^2$ and the new variable 
$u:=p^{1/\alpha'}z$, we can bound the last expression 
from below by
$$|f(0)|^2e^{-2p\psi(0)}p^{2n/\alpha'} \int_{|z|<{p^{-1/\alpha'}}}
  e^{-cp |z|^{\alpha'}}d\Leb(z)
  = |f(0)|^2e^{-2p\psi(0)}\int_{|u|<1} e^{- c|u|^{\alpha'}} d\Leb(u).$$
Therefore, we deduce from \eqref{eq:n4.19} that
$ |f(0)|^2e^{-2p\psi(0)} \lesssim p^{2n/\alpha'}$.
The estimates  \eqref{eq_extremal_Bergman_function_estimate}, 
\eqref{eq_distortion_Holder_weights} and then the theorem follow.
\endproof

In the case where $K=X$ and $\mu=\mu^0$, we have the following 
lemma.

\begin{lemma} \label{lemma_L_X}
 Let $A>0$ and  $\alpha>0$  be constants.
Let $\phi$ be an $\omega_0$-p.s.h. function on $X$ whose 
$\Cc^\alpha$-norm is bounded by $A$. 
Then there exists a constant $c_{A,\alpha}>0$ depending only 
on $X,L,h_0,A$ and $\alpha$ such that for every $p>1$, we have 
\begin{align*}
0\leq  \mL_p(\mu^0,\phi)-\mL_p(X,\phi)
\leq  {c_{A,\alpha}\log p\over p}\cdot 
\end{align*}
\end{lemma}
\proof
It is enough to apply Lemma \ref{lemma_L-K}  and 
Theorem \ref{th_BM_open} for $K=X$. Note that since $\phi$ is 
$\omega_0$-p.s.h., we have $P_X\phi=\phi$. 
\endproof

\subsection{Main estimates for the volumes and energy} 
\label{subsection_Main_estimates}

We gather in this subsection the main estimates needed for
 the proofs of our main theorems. 
 
 \medskip
\noindent
{\bf Normalization.} From now on, in order to simplify the notation, 
we use the following normalization
\begin{equation}\label{eq_normalization}
 \mE_\eq(X,0)=0\quad\text{and}
 \quad 
\mL_p(\mu^0,0)=0\quad \text{for} \quad p\in\N.
 \end{equation}
Here, the function identically 0 is used as a smooth strictly 
$\omega_0$-p.s.h. weight.

\medskip

 For continuous weights $\phi_1,\phi_2$ on $X$, the following 
 quantities will play an important role in the sequel:
\begin{align}\label{eq:n4.23}
\mV_p(\phi_1,\phi_2):= \big | \big(\mL_p(\mu^0,\phi_1)
 - \mL_p(\mu^0,\phi_2)\big)- \big(  \mE_\eq(X,\phi_1)
 - \mE_\eq(X,\phi_2)\big)\big|
\end{align}
and 
\begin{align}\label{eq:n4.24}
\mW_p(\phi_1,\phi_2):= \big | \big(\mL_p(X,\phi_1) 
- \mL_p(X,\phi_2)\big)- \big(  \mE_\eq(X,\phi_1)
- \mE_\eq(X,\phi_2)\big)\big|.
\end{align}

Here are three crucial propositions. The first two results deal with 
strictly $\omega_0$-p.s.h. weights,
whereas  the last one considers  the case  with weakly 
$\omega_0$-p.s.h. weights.

\begin{proposition}\label{prop_L_to_E_1} 
Let  $\phi_1$ and $\phi_2$ be two  weights  of class $\Cc^3$ on 
$X$ such that\break 
$\max(\|\phi_1\|_3,\|\phi_2\|_3)\leq A$ for some given constant 
$A>0$. Suppose
 $\ddc \phi_1+\omega_0\geq \zeta\omega_0$ and  
$\ddc \phi_2+\omega_0\geq \zeta\omega_0$ for some $\zeta>0$.
Then, there is a constant
$c_{A,\zeta}>0$ depending only on $X,L,\omega_0,A$ and $\zeta$
 such that for all $p> 1$
$$\mV_p(\phi_1,\phi_2)\leq c_{A,\zeta}   (\log p)^{3/2 } p^{-1/2}
\qquad \text{and}\qquad  \mW_p(\phi_1,\phi_2) 
\leq c_{A,\zeta}  (\log p)^{3/2} p^{-1/2}.$$
\end{proposition}
\proof  
By Lemma \ref{lemma_L_X}, the second estimate of
 the proposition follows from the  first one. 
 So we only need to prove the first estimate. In what follows, 
 all involved constants may depend on $X,L,\omega_0,A$ and 
 $\zeta$.
Recall that $\zeta\leq 1$ because 
$\ddc \phi_j+\omega_0\geq \zeta\omega_0$ and 
$\ddc \phi_j+\omega_0$ is cohomologous to 
$\omega_0$. It is  enough to consider $p$ large enough. 
 
For $t\in[0,1]$, define
$ \phi_t:= t\phi_1+(1-t)\phi_2$. 
By Lemma \ref{lem_directional_derivative_L}, we  get 
$$\mL_p(\mu^0,\phi_1) - \mL_p(\mu^0,\phi_2)=\int_{t=0}^1 dt
\int_X (\phi_1-\phi_2)\Bc_p(\mu^0,\phi_t).$$
Since   $\ddc \phi_t+\omega_0\geq \zeta\omega_0,$ 
by Theorem \ref{tue4} applied  to  $\phi_t$, 
the right hand side of the last identity is equal to
$$\int_{t=0}^1 dt\int_X (\phi_1-\phi_2)\mu_\eq(X,\phi_t) 
+O\big((\log p)^{3/2}p^{-1/2} \big).$$
By applying Theorem \ref{thm_differentiability}, 
the double integral in the last line  is equal to
$$\int_{t=0}^1 \left. {d\over  dt}\right|_{t=0} \mE_\eq(X,\phi_t) 
=\mE_\eq(X,\phi_1)-\mE_\eq(X,\phi_2).$$
Therefore, we get
$$\mL_p(\mu^0,\phi_1) - \mL_p(\mu^0,\phi_2) 
=  \mE_\eq(X,\phi_1)- \mE_\eq(X,\phi_2) 
 +       O \big((\log p)^{3/2}p^{-1/2}  \big) ,$$
which proves the proposition.
\endproof

\begin{proposition}\label{prop_L_to_E_2} 
Let  $0<\alpha\leq 1$ and $A>0$ be constants.  
Let  $\phi_1$ and $\phi_2$ be two  weights of class 
$\Cc^{0,\alpha}$ on $X$ such that 
$\max(\|\phi_1\|_{\Cc^{0,\alpha}},\|\phi_2\|_{\Cc^{0,\alpha}})
\leq A$. Suppose
$\ddc \phi_1+\omega_0\geq \zeta\omega_0$ and  
$\ddc \phi_2+\omega_0\geq \zeta\omega_0$ for some 
$\zeta>0.$ Then, there is a constant
$c_{A,\alpha,\zeta}>0$ depending only on
 $X,L,\omega_0,A,\alpha$ and $\zeta$ such that for all $p>1$
$$\mV_p(\phi_1,\phi_2)\leq c_{A,\alpha,\zeta}   (\log p)^{\alpha/2} 
p^{-\alpha/6}  \qquad \text{and}\qquad 
\mW_p(\phi_1,\phi_2) \leq c_{A,\alpha, \zeta}  (\log p)^{\alpha/2} 
p^{-\alpha/6}.$$
\end{proposition}
\proof
As in the last proposition, we can assume that $\zeta$ is fixed with
$\zeta\leq 1$ and $p$ is large enough. 
Moreover, we only need to prove the first estimate.
The constants involved in the calculus below may depend on 
$X,L,\omega_0,A, \alpha$ and $\zeta$. 
Fix a constant $c>0$ large enough and define 
$$\epsilon:= c\big((\log p)^{3/2} p^{-1/2}\big)^{1/3}\ll 1$$
for $p$ large enough. 
By Theorem \ref{thm_Demailly} applied to $(1-\zeta)^{-1}\phi_1$
and $(1-\zeta)^{-1}\phi_2$, 
there exist two smooth weights
$\phi_{j,\epsilon}:=(1-\zeta)\big[(1-\zeta)^{-1}
\phi_j\big]_\epsilon$  for $j=1,2$ such that
\begin{enumerate}
\item[a)] $\ddc \phi_{j,\epsilon}  +\omega_0\geq \zeta\omega_0;$
\item[b)] $\|\phi_{j,\epsilon}-\phi_j\|_\infty \lesssim 
\epsilon^\alpha$;
\item[c)] $\|\phi_{j,\epsilon}\|_{\Cc^3} \lesssim
\epsilon^{\alpha-3}$.
\end{enumerate}

We deduce from \eqref{eq:n4.23}, 
Theorem \ref{thm_differentiability} and 
Lemma \ref{lem_directional_derivative_L} that  
$$|\mV_p(\phi_1,\phi_2)
-\mV_p(\phi_{1,\epsilon},\phi_{2,\epsilon})|
\lesssim \epsilon^\alpha.$$
We can apply  Theorem \ref{tue4}  to  $\phi_{j,\epsilon}$ 
and their linear combinations as
 in the proof of Proposition \ref{prop_L_to_E_1}. 
The choice of $\epsilon$ and the above properties a)-c) allow us to 
check the hypotheses of that theorem for large $p$.  
Therefore, taking into account the estimate c),  we obtain
$$\mL_p(\mu^0,\phi_{1,\epsilon}) - \mL_p(\mu^0,\phi_{2,\epsilon})
=  \mE_\eq(X,\phi_{1,\epsilon})- \mE_\eq(X,\phi_{2,\epsilon})
+ O\big((\log p)^{3/2}p^{-1/2}  \epsilon^{ \alpha-3}  \big),$$
or equivalently
$$\mV_p(\phi_{1,\epsilon},\phi_{2,\epsilon})
\lesssim (\log p)^{3/2}p^{-1/2}  \epsilon^{ \alpha-3}.$$
Thus,
$$\mV_p(\phi_1,\phi_2)
\lesssim  (\log p)^{3/2}p^{-1/2} \epsilon^{ \alpha-3}  
+ \epsilon^{\alpha}. $$
This estimate and the choice of $\epsilon$ imply the
first inequality in the proposition.
\endproof

\begin{proposition}\label{prop_L_to_E_4} 
Let $0<\alpha\leq 1$ and $A>0$ be constants. Let 
$\phi_1$ and $\phi_2$ be two $\omega_0$-p.s.h. weights of class 
$\Cc^{0,\alpha}$  on $X$ such that 
$\max(\|\phi_1\|_{\Cc^{0,\alpha}},\|\phi_2\|_{\Cc^{0,\alpha}})
\leq A$.   Then, there is a constant
$c_{A,\alpha}>0$ depending only on $X,L,\omega_0, A$ and 
$\alpha$ such that for all  $p>1$
$$ \mV_p(\phi_1,\phi_2)\leq
c_{A,\alpha}  (\log p)^{3\beta_\alpha } p^{-\beta_\alpha} 
\qquad \text{and}\qquad 
\mW_p(\phi_1,\phi_2) \leq c_{A,\alpha}   (\log p)^{3\beta_\alpha}
p^{-\beta_\alpha},$$
where $\beta_\alpha:= \alpha/(6+3\alpha)$.
\end{proposition}
\proof  
As above, we only need to prove the first inequality and to 
consider $p$ large enough. Choose 
$$\epsilon:=(\log p)^{1/(2+\alpha)}p^{-1/(6+3\alpha)} \qquad 
\text{and} \qquad \zeta:=\epsilon^\alpha.$$
Define $\phi_j' := (1-\zeta)\phi_j$.  
We proceed as in Proposition \ref{prop_L_to_E_2} but should take 
into account the fact that $\zeta$ is no more fixed. 
The constants involved in the computation below should 
be independent of $\zeta$. 

As in that proposition, we obtain
$$|\mV_p(\phi_1,\phi_2)-\mV_p(\phi_1',\phi_2')|\lesssim \zeta$$
and since  $\ddc\phi_j' +\omega_0 \geq  \zeta\omega_0$
$$\mV_p(\phi_1',\phi_2') 
\lesssim  \zeta^{-3/2} (\log p)^{3/2}p^{-1/2} \epsilon^{\alpha-3} 
+ \epsilon^{\alpha}.$$
We then deduce that 
$$\mV_p(\phi_1,\phi_2)
\lesssim \zeta + \zeta^{-3/2}(\log p)^{3/2}p^{-1/2} 
\epsilon^{\alpha-3}  + \epsilon^{\alpha}.$$
The above choice of $\epsilon$ and $\zeta$ implies the result.
\endproof

In the rest of this subsection, we give some results which relate 
Fekete points with the functionals considered above.
Fix an orthonormal basis $S_p=(s_1,\ldots,s_{N_p})$ of 
$H^0(X,L^{p})$ with respect to the scalar product on 
$H^0(X,L^{p})$ induced by $h_0$ and $\mu^0$.
Consider a weighted compact set $(K,\phi)$ with $\phi$ 
continuous on $K$. Recall that 
$$\| \det S_p\|_{L^\infty(K,p\phi)}
:= \sup_{(x_1,\ldots,x_{N_p})\in K^{N_p}} 
|\det (s_i(x_j))|e^{-p\phi(x_1)-\cdots-p\phi(x_{N_p}) } $$
and
$$\| \det S_p\|^2_{L^2(\mu,p\phi)}
:= \int_{(x_1,\ldots,x_{N_p})\in K^{N_p}}
 |\det (s_i(x_j))|^2e^{-2p\phi(x_1)-\cdots-2p\phi(x_{N_p}) }
d\mu(x_1)\ldots d\mu(x_{N_p}),$$
if $\phi$ is a weight on $K$ and $\mu$ is a probability measure 
supported by $K$.
 
We assume further that $(K,\phi)$ is regular, i.e., 
$\phi_K=P_K\phi$, that $P_K\phi$ is continuous, and also that 
the following strong Bernstein-Markov inequality holds
\begin{equation}\label{eq_condi}
 \sup_X\rho_p(\mu^0,P_K\phi)\leq  Bp^B \quad 
 \text{for some constant}\quad  B>0.
\end{equation}

\begin{lemma}\label{lem_L_infinity_vs_L2} 
Let $S_p,K$ and $\phi$ be as above with condition 
\eqref{eq_condi}. Then there is a constant $c>0$ depending 
only on $B$ such that for $p>1$
$$\big|\log \| \det S_p\|_{L^\infty(K,p\phi)} -  
\log \| \det S_p\|_{L^2(\mu^0,pP_K\phi)} \big| \leq cN_p\log  p.$$
\end{lemma}
\proof 
Observe that the restriction of $(L^p)^{\boxtimes N_p}$ to 
$\{x_1\}\times \cdots \times\{x_{N_p-1}\}\times X$ can be
identified to the line bundle $L^p$ over $X$. Therefore, we can 
apply Proposition \ref{prop_max_principle} to 
$x\mapsto\det S_p(x_1,\ldots,x_{N_p-1},x)$. 
Then, using inductively the same argument for the other variables 
$x_i$, we obtain 
$$\| \det S_p\|_{L^\infty(K,p\phi)}
=\| \det S_p\|_{L^\infty(X,pP_K\phi)}.$$
Hence, 
$$\| \det S_p\|_{L^\infty(K,p\phi)}
\geq \| \det S_p\|_{L^2(\mu^0,pP_K\phi)}.$$
Now, to complete the proof we only need to show that
\begin{equation}\label{eq_reduction}
\log \| \det S_p\|_{L^\infty(X,pP_K\phi)}
\leq  \log  \| \det S_p\|_{L^2(\mu^0,pP_K\phi)} 
+O(N_p\log  p ).
\end{equation}
 By  \eqref{eq_condi}, we get 
$$|s(x)|^2_{pP_K\phi} \leq  \rho_p(\mu^0,P_K\phi)(x) 
\| s\|^2_{L^2(\mu^0,pP_K\phi)}
    \leq B p^B \| s\|^2_{L^2(\mu^0,pP_K\phi)}$$
for  every section $s\in H^0(X,L^{p})$, $p>1$, and  
$x\in X.$ Now, if $x_1,\ldots, x_{N_p}$ are points in $X,$ 
then for each $j$
$$x\mapsto \det S_p(x_1,\ldots,x_{j-1}, x, x_{j+1},
\ldots, x_{N_p})$$
is a holomorphic section in $H^0(X,L^p).$ 
A successive application of the last inequality for 
$j=1,2,\ldots, N_p$ yields
$$\| \det S_p\|^2_{L^\infty(X,pP_K\phi)}\leq  B^{N_p} 
p^{BN_p} \| \det S_p\|^2_{L^2(\mu^0,pP_K\phi)},$$
and  \eqref{eq_reduction} follows.
\endproof

Taking the  normalization \eqref{eq_normalization} into  account,  
 we set, for each $p>1,$
\begin{equation}\label{eq_epsilon_p}
\epsilon_p:=\big |\mL_p(\mu^0,P_K\phi)
-\mE_\eq(K,\phi)\big|= \mV_p(P_K\phi,0),
 \end{equation}
and
$$\mD_p(K,\phi):= {1\over  pN_p}
\log \| \det S_p\|_{L^\infty(K,p\phi)}.$$

\begin{proposition} \label{prop_D_L_E} 
Let $S_p, K, \phi, \epsilon_p$ and  $\mD_p(K,\phi)$ be  as  
above with condition \eqref{eq_condi}.  Then there is a constant 
$c>0$ depending only on $X,L$ and $B$ such that for $p>1$
$$|\mD_p(K,\phi)+\mE_\eq(K,\phi)|
\leq c  (p^{-1}\log p+\epsilon_p),$$
 and for any Fekete measure $\mu_p$ associated with $(K,\phi)$
$$|\mL_p(\mu_p,\phi)-\mE_\eq(K,\phi)|\leq c (p^{-1}\log p
+\epsilon_p).$$
\end{proposition}
\proof 
We prove the first assertion. 
By  Lemma \ref{lem_L_infinity_vs_L2}, we only need to check that 
\begin{equation}\label{eq_reduction_L2_vs_Linfty}
 \Big|{1\over pN_p} \log \| \det S_p\|_{L^2(\mu^0,pP_K\phi)} 
 +\mE_\eq(K,\phi)\Big| \lesssim  p^{-1}\log p+\epsilon_p.
\end{equation}
Using that $S_p$ is an orthonormal basis, a direct computation 
(see \cite [Lemma 5.3]{BermanBoucksom} and 
\cite[p.377]{BermanBoucksom}), gives
$$\| \det S_p\|^2_{L^2(\mu^0,pP_K\phi)}
= N_p! {  \vol \mathcal B^2_p(\mu^0,0)\over
  \vol \mathcal B^2_p(\mu^0,P_K\phi)},$$
which  implies
$${1\over pN_p} \log \| \det S_p\|_{L^2(\mu^0,pP_K\phi)} 
= \mL_p(\mu^0,0)- \mL_p(\mu^0,P_K\phi) 
+ {\log  N_p!\over 2pN_p}\cdot$$
By the normalization \eqref{eq_normalization}  
and  \eqref{eq_epsilon_p},
$$\mL_p(\mu^0,0)=0  \quad \text{and}\quad   
\mL_p(\mu^0,P_K\phi)=  \mE_\eq(K,\phi) \pm \epsilon_p.$$
 On the other hand, since  $N_p\simeq p^n$ by \eqref{e:Demailly}, 
 we have
$${\log N_p!\over 2pN_p}\lesssim {p^n\log p\over 2pN_p}
\lesssim p^{-1}\log p.$$
 Combining the last four estimates together, we obtain
  \eqref{eq_reduction_L2_vs_Linfty}.

Consider now the second assertion in the proposition. 
Using the definition of Fekete points, we obtain 
(see  \cite[(2.4)]{BBW})
$${1\over  2p N_p}\log{  \vol \mathcal B^2_p(\mu^0,0) 
\over \vol\mathcal B^2_p(\mu_p,\phi)}
=\mD_p(K,\phi)-{1\over 2p} \log N_p.$$
By the normalization \eqref{eq_normalization}, the left-hand side is 
 $-\mL_p(\mu_p,\phi).$ Using again that $N_p\simeq p^n, $  
 we deduce the result from the first assertion of the proposition.
\endproof

\subsection{Proofs of the main results and further remarks}

In this subsection, we will give the proofs of the main theorems 
stated in the Introduction. We need the following auxiliary lemmas.

\begin{lemma}\label{lem_MA_Lip} There is  a constant $c>0$ such 
that for every  continuous weight  $\phi$ on $K$  and every 
function $v$ of class $\Cc^{1,1}$ on $X$, we have   
$$\big|\langle \mu_\eq(K, \phi+tv) -\mu_\eq(K, \phi),v\rangle\big|
\leq  c|t|\|v\|_{L^\infty(K)}\|\ddc v\|_\infty \quad \text{for}
\quad t\in\R.$$
\end{lemma}
\proof 
Define 
$$ \Psi:=\sum_{j=1}^n (\ddc P_K\phi+\omega_0)^{j-1}
\wedge(\ddc P_K(\phi+tv)+\omega_0)^{n-j}.$$
Observe that $\ddc P_K\phi+\omega_0$ and  
$\ddc P_K(\phi+tv)+\omega_0$ 
are positive closed $(1,1)$-currents cohomologous to 
$\omega_0$. So $\Psi$ is a sum of $n$ positive closed 
$(n-1,n-1)$-currents of bounded mass. 
Define also $u:=P_K(\phi+tv)-P_K\phi$.
For  $t\in\R,$ we have
\begin{eqnarray*}
\langle \mu_\eq(K,\phi+tv),v\rangle  
-  \langle \mu_\eq(K,\phi),v\rangle  
& = & \big \langle \NMA(P_K(\phi+tv))-\NMA(P_K\phi), v\big\rangle  \\
& = & \const \langle \ddc u\wedge \Psi,v\rangle \ 
= \const   \langle\ddc v\wedge\Psi, u  \rangle.
\end{eqnarray*}
On the  other hand, 
by Lemma \ref{lem_projection_P_K_is_Lipschitz},
$$\|u\|_{L^\infty(X)}=\|P_K(\phi+tv)-P_K\phi\|_{L^\infty(X)}
\leq  |t|\|v\|_{L^\infty(K)}.$$
Since $v\in\Cc^{1,1}(X),$ $\ddc v$ can be written as the 
difference of two positive closed bounded $(1,1)$-forms.
Consequently, $\ddc v\wedge\Phi$ is a signed sum of $2n$ 
positive measures of bounded mass. This and the above 
computation imply the lemma.
 \endproof

\begin{lemma}\label{key_lemma}
Let $ \epsilon>0$ and $M>0$ be constants. Let  $F$ and $G$ be   
functions defined on $[-\epsilon^{1/2},\epsilon^{1/2}]$  such that
\begin{enumerate}
\item[a)] $F(t)\geq G(t)-\epsilon$ and 
$|F(0) -  G(0)|\leq \epsilon;$
\item[b)]  $F$ is  concave on $[-\epsilon^{1/2},\epsilon^{1/2}]$ 
and differentiable at $0;$
\item[c)] $G$ is differentiable in $ [-\epsilon^{1/2},\epsilon^{1/2}],$ 
and its derivative $G'$ satisfies
$|G'(t)-G'(0)|\leq M\epsilon^{1/2}$ for $t \in [-\epsilon^{1/2},
  \epsilon^{1/2}] .$ The last inequality holds when 
  $|G'(t)-G'(0)|\leq M|t|$. 
\end{enumerate}
Then we have
$$|F'(0)- G'(0)|\leq (2+M)\epsilon^{1/2}.$$
\end{lemma}

\proof 
This is a quantitative version of \cite[Lemma 7.6]{BermanBoucksom}.
Since $F$ is concave, we have
$$F(0)+F'(0)t\geq  F(t)$$
for $|t|\leq \epsilon^{1/2}.$  Hence, for $t:=\pm\epsilon^{1/2},$ 
we get   
\begin{equation}\label{eq_key_lemma}
 t F'(0)\geq  G(t)-G(0)-2\epsilon=G(t)-G(0)-2t^2.
\end{equation}
Now, take $t:=\epsilon^{1/2}.$ There exists $s\in (0,t)$ such that
$${G(t)-G(0)\over t}= G'(s) \qquad\text{and by c)}
 \qquad  |G'(s)-G'(0)|\leq Mt.$$
This, combined with \eqref{eq_key_lemma} yields
$$ F'(0)\geq  G'(s)-  2t \geq   G'(0)  - (2+M)t.$$
Hence, 
$F'(0)-  G'(0)  \geq  -  (2+M)\epsilon^{1/2}.$
The inequality  
$ F'(0)-  G'(0)  \leq   (2+M)\epsilon^{1/2}$ is obtained in the same 
way by using $t:=-\epsilon^{1/2}$.
\endproof

\medskip
 \noindent {\bf End of the proof of  
 Theorem \ref{thm_main_X_positive}.}
By \eqref{eq:n1.23}, we only need to consider the case 
$\gamma=3,$ i.e., to prove
 \begin{equation}\label{eq_main_theorem_alpha_2}
   \big |  \langle \mu_p -\mu_\eq(X,\phi), v \rangle \big | \lesssim    
   p^{-1/4 }(\log p)^{3/4} 
 \end{equation}
 for every test function $v$ such that
 $\| v\|_{\Cc^{3}}\leq 1.$
We will apply Lemma \ref{key_lemma} to the following functions
$$F(t):=\mL_p(\mu_p,\phi+tv)\qquad \text{and}
\qquad  G(t):= \mE_\eq(X,\phi+tv).$$

By  Lemma \ref{lemma_L_K_mu},   
\begin{align}\label{eq:n4.89}
\mL_p(\mu_p,\phi+tv)\geq  \mL_p(X,\phi+tv).
\end{align}
On the other hand, since $\ddc  v$ is bounded, we can find 
a constant $t_0>0$ such that $\phi+tv$ is 
$(1-\zeta)\omega_0$-p.s.h. for $|t|\leq t_0$ and $\zeta>0$ 
a fixed constant. Recall  that the function $0$ satisfies 
the normalization \eqref{eq_normalization}. Consequently,
Proposition \ref {prop_L_to_E_1},  applied to $\phi+tv$ and the 
function $0,$ yields 
$$ |\mL_p(X,\phi+tv)-\mE_\eq(X,\phi+tv)|
\lesssim  p^{-1/2 }(\log p)^{3/2}.$$
This, combined with \eqref{eq:n4.89},
 shows that
 \begin{align}\label{eq:n4.91}
 F(t)- G(t) \gtrsim -   p^{-1/2 }(\log p)^{3/2}  . 
\end{align}

Next, since $\phi$ is $\omega_0$-p.s.h., we have $P_K\phi=\phi$. 
Moreover, 
we have the strong Bernstein-Markov inequality thanks to
Theorem \ref{th_BM_open} applied to $K:=X$.
Let $\epsilon_p$ be defined as
in \eqref{eq_epsilon_p} with $K=X$ and $P_K\phi=\phi$.
By Proposition \ref{prop_L_to_E_1} again, we have 
$\epsilon_p= O(p^{-1/2 }(\log p)^{3/2})$.
Consequently, applying Proposition \ref{prop_D_L_E} yields 
\begin{align}\label{eq:n4.92} 
 |F(0)-G(0)| \lesssim p^{-1/2 }(\log p)^{3/2}.
\end{align}
Recall from Lemma \ref{lem_L_concave} that $F$ is concave.
Moreover, by Lemma \ref{lem_directional_derivative_L}, we have
\begin{align}\label{eq:n4.93}
F'(0)=\langle v, \Bc_p(\mu_p,\phi)  \rangle.
\end{align}
On the other hand, by Theorem \ref{thm_differentiability}, $G$
 is differentiable with
\begin{align}\label{eq:n4.94}
G'(t)=\langle v, \mu_\eq( X,\phi+tv)  \rangle.
\end{align}

Finally, by  Lemma \ref{lem_MA_Lip}, condition c) in  
Lemma \ref{key_lemma} is satisfied for a suitable constant $M>0$. 
Combining this and the discussion between 
\eqref{eq:n4.91}-\eqref{eq:n4.94}, we are in the position to apply 
Lemma \ref{key_lemma} to a constant $\epsilon$ of order 
$p^{-1/2}(\log p)^{3/2}$.
Using  the  above expression for $F'(0)$ and $G'(0),$  we get 
$$\big | \langle \Bc_p(\mu_p,\phi) , v \rangle
 - \langle \mu_\eq(X, \phi), v  \rangle \big | 
 =O  \big(p^{-1/4} (\log p)^{3/4}  \big).$$
Recall from the discussion before 
Lemma \ref{lem_directional_derivative_L} that 
$\Bc_p(\mu_p,\phi)= \mu_p.$ Hence,  estimate  
\eqref{eq_main_theorem_alpha_2} follows immediately.
 \hfill $\square$

\begin{remark} \rm \label{rk_main_X_positive}
If in Theorem \ref{thm_main_X_positive}, the function $\phi$ is only 
$\Cc^{0,\alpha}$ for some $0<\alpha\leq 1$, we can apply 
Proposition \ref{prop_L_to_E_2} instead of \ref{prop_L_to_E_1}
in order to get 
$$\dist_\gamma(\mu_p,\mu_\eq(X,\phi))
\lesssim (\log p)^{\alpha\gamma/8}p^{-\alpha\gamma/24} 
\quad \text{for} \quad 0<\gamma\leq 2.$$
\end{remark}

\medskip
 \noindent {\bf End of the proof of Theorem \ref{thm_main}.}
By \eqref{eq:n1.23}, we only need to consider the case 
$\gamma=2,$ i.e., to prove
 \begin{equation*}
   \big |  \langle \mu_p -\mu_\eq(K,\phi), v \rangle \big | 
   \lesssim    p^{-2\beta} (\log p)^{6\beta}  
 \end{equation*}
 for every test $\Cc^2$ function $v$ such that
 $\| v\|_{\Cc^2}\leq 1.$ Recall that $\beta:=\alpha'/(24+12\alpha')$.
Define 
$$F(t):=\mL_p(\mu_p,\phi+tv)\qquad \text{and}\qquad  G(t)
:= \mE_\eq(K,\phi+tv)=\mE_\eq(X,P_K(\phi+tv))$$
for $t$ in a neighborhood of $0\in\R$. 
By Lemma \ref{lemma_L_K_mu} and 
Proposition \ref{prop_max_principle},
 $$\mL_p(\mu_p,\phi+tv)\geq  \mL_p(K,\phi+tv)
 = \mL_p(X,P_K(\phi+tv)).$$
As $0<\alpha\leq 2,$ we infer that
$\phi+tv\in\Cc^\alpha(K).$
Since $K$ is $(\Cc^\alpha,\Cc^{\alpha'})$-regular,  we deduce 
that $P_K(\phi+tv)$ is an 
$\omega_0$-p.s.h. weight on $X$ with bounded 
$\Cc^{\alpha'}$-norm.
This, coupled with 
Proposition \ref{prop_L_to_E_4}, applied to $P_K(\phi+tv)$ 
and the function 0,  for $\alpha'$ instead of $\alpha$, and 
the normalization \eqref{eq_normalization}, shows that
$$  F(t)- G(t) \gtrsim -   p^{-4\beta} (\log p)^{12\beta}.$$

By Theorem \ref{th_BM_open},
condition \eqref{eq_condi} is fulfilled.
Let $\epsilon_p$ be defined as in \eqref{eq_epsilon_p}.
By  Proposition \ref{prop_L_to_E_4} for $\alpha'$ instead of 
$\alpha$, we have $\epsilon_p=O(p^{-4\beta} (\log p)^{12\beta})$. 
Consequently, applying Proposition \ref{prop_D_L_E} yields   
$$|F(0)-G(0)| \lesssim    p^{-4\beta} (\log p)^{12\beta}.$$

Finally, since  $\|v\|_{\Cc^2(X)}\leq 1,$ we can check condition c) 
in Lemma \ref{key_lemma} using 
  Lemma \ref{lem_MA_Lip}. Applying Lemma \ref{key_lemma} 
  to a constant $\epsilon$ of order 
  $p^{-4\beta} (\log p)^{12\beta}$, we easily obtain the result 
  as in  the proof of Theorem \ref{thm_main_X_positive}.
 \hfill $\square$

\begin{remark} \rm
Optimal estimates for the speed of convergence in our results are
still unknown. This is an interesting problem which may require 
a better understanding of the Bergman kernels. Results in this 
direction may have consequences in theory of sampling and 
interpolation for line bundles with singular metric and not necessarily 
of positive curvature. Demailly suggested us to study first the case
in $\C^n$ with data invariant under the action of the real torus 
$(\S^1)^n$. 
\end{remark}

\begin{remark} \rm
Our proofs still hold for {\it almost} Fekete configurations 
$P=(x_1,\ldots,x_{N_p})\in K^{N_p}$ in the sense that the quantity
$\sigma_P$ below is not too big. Assume that $P$ is not 
necessarily a Fekete configuration and define 
$$\sigma_P:={1\over pN_p} \log \|\det S_p\|_{L^\infty(K,p\phi)}
-{1\over pN_p} \log \|\det S_p(P)\|_{p\phi}.$$
Then our main estimates are still valid for this configuration
if we add to their right hand sides the term 
$O(\sigma_P^{\gamma/4})$ for the estimates in 
Theorems \ref{th_main_Cn}, \ref{thm_main} and 
Corollary \ref{cor_main_K}, and $O(\sigma_P^{\gamma/6})$ 
for the estimate in Theorem \ref{thm_main_X_positive}. 
The main change in the proofs is that we need to add  
$O(\sigma_P)$ to the right hand side of the second inequality
in Proposition \ref{prop_D_L_E}. This answers a question that
Norm Levenberg asked us, see also \cite{Levenberg2}. 
\end{remark}

\end{document}